\documentclass[12pt]{article}

\usepackage{amsmath,amssymb,amsthm,amscd,a4wide,upref}

\usepackage[enableskew]{youngtab}
\usepackage{ytableau}
\usepackage{tikz}
\usepackage{mathtools}

\usepackage{hyperref}
\hypersetup{colorlinks,citecolor=blue,filecolor=black,linkcolor=blue,urlcolor=blue}

\usepackage{lmodern}     
\usepackage[T1]{fontenc}

\begin{document}


\newcommand{\non}{\nonumber}
\newcommand{\scl}{\scriptstyle}
\newcommand{\sclnearrow}{{\scl\nearrow}\ts}
\newcommand{\scloplus}{{\scl\bigoplus}}
\newcommand{\wt}{\widetilde}
\newcommand{\wh}{\widehat}
\newcommand{\ot}{\otimes}
\newcommand{\fand}{\quad\text{and}\quad}
\newcommand{\Fand}{\qquad\text{and}\qquad}
\newcommand{\ts}{\,}
\newcommand{\tss}{\hspace{1pt}}
\newcommand{\lan}{\langle\ts}
\newcommand{\ran}{\ts\rangle}
\newcommand{\vl}{\tss|\tss}
\newcommand{\qin}{q^{-1}}
\newcommand{\tpr}{t^{\tss\prime}}
\newcommand{\spr}{s^{\tss\prime}}
\newcommand{\di}{\partial}
\newcommand{\hra}{\hookrightarrow}
\newcommand{\antiddots}
    {\underset{\displaystyle\cdot\quad\ }
    {\overset{\displaystyle\quad\ \cdot}{\cdot}}}
\newcommand{\dddots}
    {\underset{\displaystyle\quad\ \cdot}
    {\overset{\displaystyle\cdot\quad\ }{\cdot}}}
\newcommand{\atopn}[2]{\genfrac{}{}{0pt}{}{#1}{#2}}

\newcommand{\su}{s^{}}
\newcommand{\vac}{\mathbf{1}}
\newcommand{\vacf}{\tss|0\rangle}
\newcommand{\BL}{ {\overline L}}
\newcommand{\BD}{ {\overline D}}
\newcommand{\BE}{ {\overline E}}
\newcommand{\BP}{ {\overline P}}
\newcommand{\ol}{\overline}
\newcommand{\pr}{^{\tss\prime}}
\newcommand{\ba}{\bar{a}}
\newcommand{\bb}{\bar{b}}
\newcommand{\eb}{\bar{e}}
\newcommand{\bi}{\bar{\imath}}
\newcommand{\bj}{\bar{\jmath}}
\newcommand{\bh}{\bar{h}}
\newcommand{\bk}{\bar{k}}
\newcommand{\bl}{\bar{l}}
\newcommand{\hb}{\mathbf{h}}
\newcommand{\gb}{\mathbf{g}}
\newcommand{\For}{\qquad\text{or}\qquad}
\newcommand{\OR}{\qquad\text{or}\qquad}
\newcommand{\emp}{\mbox{}}


\newcommand{\U}{{\rm U}}
\newcommand{\Z}{{\rm Z}}
\newcommand{\ZY}{{\rm ZY}}
\newcommand{\Ar}{{\rm A}}
\newcommand{\Br}{{\rm B}}
\newcommand{\Cr}{{\rm C}}
\newcommand{\Fr}{{\rm F}}
\newcommand{\hF}{{\wh F}}
\newcommand{\hE}{{\wh E}}
\newcommand{\Mr}{{\rm M}}
\newcommand{\Sr}{{\rm S}}
\newcommand{\Prm}{{\rm P}}
\newcommand{\Lr}{{\rm L}}
\newcommand{\Ir}{{\rm I}}
\newcommand{\Jr}{{\rm J}}
\newcommand{\Qr}{{\rm Q}}
\newcommand{\Rr}{{\rm R}}
\newcommand{\X}{{\rm X}}
\newcommand{\Y}{{\rm Y}}
\newcommand{\DY}{ {\rm DY}}
\newcommand{\Or}{{\rm O}}
\newcommand{\SO}{{\rm SO}}
\newcommand{\GL}{{\rm GL}}
\newcommand{\Spr}{{\rm Sp}}
\newcommand{\Zr}{{\rm Z}}
\newcommand{\ev}{{\rm ev}}
\newcommand{\op}{{\rm op}}
\newcommand{\Pf}{{\rm Pf}}
\newcommand{\Ann}{{\rm{Ann}\ts}}
\newcommand{\Norm}{{\rm Norm\tss}}
\newcommand{\Ad}{{\rm Ad}}
\newcommand{\SY}{{\rm SY}}
\newcommand{\Pff}{{\rm Pf}\tss}
\newcommand{\Hf}{{\rm Hf}\tss}
\newcommand{\trts}{{\rm tr}\ts}
\newcommand{\aff}{{\rm aff}}
\newcommand{\otr}{{\rm otr}}
\newcommand{\cri}{{\rm cri}}
\newcommand{\row}{{\rm row}}
\newcommand{\End}{{\rm{End}\ts}}
\newcommand{\Mat}{{\rm{Mat}}}
\newcommand{\Hom}{{\rm{Hom}}}
\newcommand{\id}{{\rm id}}
\newcommand{\middd}{{\rm mid}}
\newcommand{\ch}{{\rm{ch}\ts}}
\newcommand{\ind}{{\rm{ind}\ts}}
\newcommand{\Normts}{{\rm{Norm}\ts}}
\newcommand{\mult}{{\rm{mult}}}
\newcommand{\per}{{\rm per}\ts}
\newcommand{\sgn}{{\rm sgn}\ts}
\newcommand{\sign}{{\rm sign}\ts}
\newcommand{\qdet}{{\rm qdet}\ts}
\newcommand{\sdet}{{\rm sdet}\ts}
\newcommand{\Ber}{{\rm Ber}\ts}
\newcommand{\inv}{{\rm inv}\ts}
\newcommand{\inva}{{\rm inv}}
\newcommand{\grts}{{\rm gr}\ts}
\newcommand{\grpr}{{\rm gr}^{\tss\prime}\ts}
\newcommand{\degpr}{{\rm deg}^{\tss\prime}\tss}
\newcommand{\Cond}{ {\rm Cond}\tss}
\newcommand{\Fun}{{\rm{Fun}\ts}}
\newcommand{\Rep}{{\rm{Rep}\ts}}
\newcommand{\sh}{{\rm{sh}}}
\newcommand{\weight}{{\rm{wt}\ts}}
\newcommand{\chara}{{\rm{char}\ts}}
\newcommand{\diag}{ {\rm diag}}
\newcommand{\Bos}{ {\rm Bos}}
\newcommand{\Ferm}{ {\rm Ferm}}
\newcommand{\cdet}{ {\rm cdet}}
\newcommand{\rdet}{ {\rm rdet}}
\newcommand{\imm}{ {\rm imm}}
\newcommand{\ad}{ {\rm ad}}
\newcommand{\tr}{ {\rm tr}}
\newcommand{\gr}{ {\rm gr}\tss}
\newcommand{\str}{ {\rm str}}
\newcommand{\loc}{{\rm loc}}
\newcommand{\Gr}{{\rm G}}

\newcommand{\twobar}{{\bar 2}}
\newcommand{\threebar}{{\bar 3}}


\newcommand{\AAb}{\mathbb{A}\tss}
\newcommand{\CC}{\mathbb{C}}
\newcommand{\FF}{\mathbb{F}}
\newcommand{\KK}{\mathbb{K}\tss}
\newcommand{\QQ}{\mathbb{Q}\tss}
\newcommand{\SSb}{\mathbb{S}\tss}
\newcommand{\TT}{\mathbb{T}\tss}
\newcommand{\ZZ}{\mathbb{Z}\tss}
\newcommand{\Sbb}{\mathbb{S}}
\newcommand{\ZZb}{\mathbb{Z}}


\newcommand{\Ac}{{\mathcal A}}
\newcommand{\Bc}{{\mathcal B}}
\newcommand{\Cc}{{\mathcal C}}
\newcommand{\Cl}{{\mathcal Cl}}
\newcommand{\Dc}{{\mathcal D}}
\newcommand{\Ec}{{\mathcal E}}
\newcommand{\Fc}{{\mathcal F}}
\newcommand{\Jc}{{\mathcal J}}
\newcommand{\Gc}{{\mathcal G}}
\newcommand{\Hc}{{\mathcal H}}
\newcommand{\Lc}{{\mathcal L}}
\newcommand{\Nc}{{\mathcal N}}
\newcommand{\Xc}{{\mathcal X}}
\newcommand{\Yc}{{\mathcal Y}}
\newcommand{\Oc}{{\mathcal O}}
\newcommand{\Pc}{{\mathcal P}}
\newcommand{\PD}{{\mathcal {PD}}}
\newcommand{\Qc}{{\mathcal Q}}
\newcommand{\Rc}{{\mathcal R}}
\newcommand{\Sc}{{\mathcal S}}
\newcommand{\Tc}{{\mathcal T}}
\newcommand{\Uc}{{\mathcal U}}
\newcommand{\Vc}{{\mathcal V}}
\newcommand{\Wc}{{\mathcal W}}
\newcommand{\Zc}{{\mathcal Z}}
\newcommand{\HC}{{\mathcal HC}}


\newcommand{\asf}{\mathsf a}
\newcommand{\bsf}{\mathsf b}
\newcommand{\csf}{\mathsf c}
\newcommand{\nsf}{\mathsf n}


\newcommand{\Sym}{\mathfrak S}
\newcommand{\h}{\mathfrak h}
\newcommand{\q}{\mathfrak q}
\newcommand{\n}{\mathfrak n}
\newcommand{\m}{\mathfrak m}
\newcommand{\p}{\mathfrak p}
\newcommand{\gl}{\mathfrak{gl}}
\newcommand{\oa}{\mathfrak{o}}
\newcommand{\spa}{\mathfrak{sp}}
\newcommand{\osp}{\mathfrak{osp}}
\newcommand{\g}{\mathfrak{g}}
\newcommand{\kgot}{\mathfrak{k}}
\newcommand{\agot}{\mathfrak{a}}
\newcommand{\bgot}{\mathfrak{b}}
\newcommand{\sll}{\mathfrak{sl}}
\newcommand{\f}{\mathfrak{f}}
\newcommand{\z}{\mathfrak{z}}
\newcommand{\Zgot}{\mathfrak{Z}}


\newcommand{\al}{\alpha}
\newcommand{\be}{\beta}
\newcommand{\ga}{\gamma}
\newcommand{\de}{\delta}
\newcommand{\De}{\Delta}
\newcommand{\Ga}{\Gamma}
\newcommand{\ep}{\epsilon}
\newcommand{\ee}{\epsilon^{}}
\newcommand{\ve}{\varepsilon}
\newcommand{\ls}{\ts\lambda\ts}
\newcommand{\vk}{\varkappa}
\newcommand{\vs}{\varsigma}
\newcommand{\vt}{\vartheta}
\newcommand{\ka}{\kappa}
\newcommand{\vp}{\varphi}
\newcommand{\la}{\lambda}
\newcommand{\La}{\Lambda}
\newcommand{\si}{\sigma}
\newcommand{\ta}{\theta}
\newcommand{\ze}{\zeta}
\newcommand{\om}{\omega}
\newcommand{\Om}{\Omega}
\newcommand{\up}{\upsilon}


\newtheorem{thm}{Theorem}[section]
\newtheorem{lemma}[thm]{Lemma}
\newtheorem{prop}[thm]{Proposition}
\newtheorem{cor}[thm]{Corollary}
\newtheorem{conj}[thm]{Conjecture}

\theoremstyle{definition}
\newtheorem{definition}[thm]{Definition}
\newtheorem{example}[thm]{Example}

\theoremstyle{remark}
\newtheorem{remark}[thm]{Remark}

\newcommand{\bth}{\begin{thm}}
\renewcommand{\eth}{\end{thm}}
\newcommand{\bpr}{\begin{prop}}
\newcommand{\epr}{\end{prop}}
\newcommand{\ble}{\begin{lemma}}
\newcommand{\ele}{\end{lemma}}
\newcommand{\bco}{\begin{cor}}
\newcommand{\eco}{\end{cor}}
\newcommand{\bex}{\begin{example}}
\newcommand{\eex}{\end{example}}
\newcommand{\bde}{\begin{definition}}
\newcommand{\ede}{\end{definition}}
\newcommand{\bre}{\begin{remark}}
\newcommand{\ere}{\end{remark}}
\newcommand{\bcj}{\begin{conj}}
\newcommand{\ecj}{\end{conj}}

\renewcommand{\theequation}{\arabic{section}.\arabic{equation}}

\numberwithin{equation}{section}


\newcommand{\bpf}{\begin{proof}}
\newcommand{\epf}{\end{proof}}


\def\beql#1{\begin{equation}\label{#1}}

\newcommand{\bal}{\begin{aligned}}
\newcommand{\eal}{\end{aligned}}
\newcommand{\beq}{\begin{equation}}
\newcommand{\eeq}{\end{equation}}
\newcommand{\ben}{\begin{equation*}}
\newcommand{\een}{\end{equation*}}

\title{\Large\bf Harish-Chandra images of orthosymplectic Sugawara operators and Casimir elements}

\author{{Alexander Molev, \quad Madeline Nurcombe \quad and\quad Eric Ragoucy}}

\date{} 
\maketitle


\begin{abstract}
We consider the recently constructed Segal--Sugawara vectors
for the orthosymplectic Lie superalgebras. We calculate their images with respect to the Harish-Chandra
homomorphism and extend this calculation to the associated Sugawara operators and Casimir elements.
We also produce higher Gaudin Hamiltonians and elements of the quantum shift-of-argument subalgebras
in the orthosymplectic enveloping algebra. In the Appendix, we review analogous results for
the general linear Lie superalgebras.



\end{abstract}

\section{Introduction}
\label{sec:int}

In a recent work by two of us \cite{mn:ss}, a family of {\em Segal--Sugawara vectors} $\Phi_2, \Phi_3,\dots$ associated
with the orthosymplectic Lie superalgebra $\osp_{M|2n}$ was constructed. The $\Phi_k$ are elements
of the centre of the affine vertex algebra at the critical level associated with $\osp_{M|2n}$.
More specifically, setting $\g=\osp_{M|2n}$, consider
the corresponding
affine Kac--Moody superalgebra $\wh\g$ defined as a central extension
\beql{km}
\wh\g=\g\tss[t,t^{-1}]\oplus\CC K
\eeq
of the Lie
superalgebra of Laurent polynomials $\g\tss[t,t^{-1}]$.
The {\em vacuum module $V_{\cri}(\g)$ at the critical level}
is the quotient of the
universal enveloping algebra $\U(\wh\g)$ by the left
ideal generated by $\g[t]$ and $K+h^\vee$,
where $h^\vee=M-2n-2$.
The vacuum module has a vertex algebra structure and is known as
the {\em affine vertex algebra}; see e.g.
\cite{f:lc} and \cite{k:va} for definitions. The {\em centre}
of this vertex algebra is defined as the subspace
\beql{centsup}
\z(\wh\g)=\{S\in V_{\cri}(\g)\ |\ \g[t]\ts S=0\}.
\eeq
The centre has the structure of
a commutative associative superalgebra and it
can be regarded as a subalgebra of $\U(t^{-1}\g\tss[t^{-1}])$.
Any element of $\z(\wh\g)$
is called a {\em Segal--Sugawara vector\/}.
The superalgebra $\z(\wh\g)$ is equipped with the derivation $\tau=-d/dt$
arising from the vertex algebra structure.

In this paper, we use the Harish-Chandra homomorphism as
an instrument for understanding the centre $\z(\wh\g)$.
We calculate the images of the Segal--Sugawara vectors $\Phi_k$
under the Harish-Chandra
homomorphism as polynomials in generators of the algebra $t^{-1}\h[t^{-1}]$,
where $\h$ is a Cartan subalgebra of $\g$ (Theorem~\ref{thm:hch}).
Our methods rely on the previous work \cite{mm:yc}
devoted to the orthogonal and symplectic Lie superalgebras; see also \cite{m:so}.

It was conjectured in \cite[Conjecture~2.3]{mn:ss} that if $M$ is odd, then the family $\Phi_2,\Phi_4,\Phi_6\dots$
of Segal--Sugawara vectors generates the superalgebra $\z(\wh\g)$.
Relying on Theorem~\ref{thm:hch},
we extend this conjecture to all values of $M$ and describe the size
of $\z(\wh\g)$ by making a connection with the affine supersymmetric functions;
see Conjectures~\ref{conj:chevalleyosp}
and \ref{conj:even} below.

The affine vertex algebra associated with the orthosymplectic
Lie superalgebra $\osp_{1|2}$
was previously investigated in \cite{a:rs}.
A conjectural description of its centre at the critical level
was pointed out in \cite[Remark~10]{a:rs}.
A few families of Segal--Sugawara vectors for the general linear Lie superalgebra
$\gl_{m|n}$ were constructed in an earlier paper \cite{mr:mm}, and it was conjectured therein that
these vectors generate the centre of the associated
affine vertex algebra at the critical level. A further description of the size of the
centre in terms of plane partitions was conjectured in \cite{mm:iv} and confirmed in the
particular case of $\gl_{1|1}$. More recently, the conjectures have been confirmed
in \cite{afn:ca} for $n=1$, generalizing the previous work \cite{an:ca} for $\gl_{2|1}$.
We will review the results of \cite{mr:mm} in Appendix~\ref{sec:gl} below and
produce the Harish-Chandra images of the Segal--Sugawara vectors for
$\gl_{m|n}$ by using noncommutative supersymmetric polynomials (Theorem~\ref{thm:hcha})
and connecting them with the pseudo-differential operators used in \cite{afn:ca}.

To compare the properties of the centre in the super and non-super cases, recall that
if $\g$ is a simple Lie algebra, the structure of the corresponding centre $\z(\wh\g)$
is described by the theorem of Feigin and Frenkel~\cite{ff:ak}.
The differential algebra
$\z(\wh\g)$ is generated by elements
$S_1,\dots,S_n$ so that
$\z(\wh\g)$ is the algebra of polynomials
\ben
\z(\wh\g)=\CC[\tau^{\tss r}S_l\ |\ l=1,\dots,n,\ \ r\geqslant 0],
\een
where $n=\text{rank}\ts\g$. Furthermore,
the centre can be identified
with the {\em classical $\Wc$-algebra}
associated with the Langlands dual Lie algebra ${}^L\g$
via an affine version of the Harish-Chandra isomorphism
\beql{hchiaff}
\z(\wh\g)\cong \Wc({}^L\g);
\eeq
see \cite{f:lc}. A direct proof of the isomorphism \eqref{hchiaff}
for the general linear and orthogonal Lie algebras
can also be obtained from the explicit constructions of generators
of the Feigin--Frenkel centre $\z(\wh\g)$
and classical $\Wc$-algebras; see \cite[Sec.~13.1]{m:so}.
Here we get a similar independent proof of \eqref{hchiaff}
for the symplectic Lie algebras $\spa_{2n}$
by calculating the Harish-Chandra images of generators of $\z(\wh\spa_{2n})$
in Theorem~\ref{thm:sphch}; see Remark~\ref{rem:typec}.

Returning to the Lie superalgebras, note that a general super-version of the
isomorphism \eqref{hchiaff} is unknown. Since the centre consists of even elements,
one could expect that the Harish-Chandra image of $\z(\wh\g)$ is properly embedded into
the classical $\Wc$-superalgebra associated with ${}^L\g$; see e.g. \cite{s:ca}.

As an application of our main theorem (Theorem~\ref{thm:hch}), we produce new generators of the centre
of the universal enveloping algebra $\U(\osp_{M|2n})$ and calculate their images
as supersymmetric polynomials under the
Harish-Chandra isomorphism (Theorem~\ref{thm:hchfd});
see \cite[Sec.~2.2]{cw:dr}. To our knowledge, the earliest
construction of Casimir elements for $\osp_{M|2n}$ goes back to Scheunert~\cite{s:ec}
where Gelfand-type invariants were produced and their Harish-Chandra images
calculated; see also \cite{lw:sw} for a more detailed discussion of
the orthosymplectic invariants and Casimir elements.

By using the vertex algebra structure of the vacuum module $V_{\cri}(\g)$,
we apply the state-field correspondence map to the vectors $\Phi_k$
to produce the corresponding {\em Sugawara operators}
as elements of the centre of the completed universal enveloping algebra
$\wt\U_{\cri}(\wh\g)$ at the critical level; cf. \cite[Sec.~4.3]{f:lc} and \cite[Secs~7.2, 8.5]{m:so}.
We then use a super-version of the affine Harish-Chandra homomorphism to
calculate their images for both the general linear and orthosymplectic
Lie superalgebras (Corollaries~\ref{cor:hchss} and \ref{cor:hchssgln});
cf. \cite[Sec.~13.3]{m:so}.

Furthermore, we show that the construction of higher Gaudin Hamiltonians
originating in \cite{ffr:gm}
extends to the orthosymplectic Lie superalgebras (Proposition~\ref{prop:gh});
cf. \cite{mr:mm} and \cite[Ch.~14]{m:so}. Then
we follow the approach of \cite{r:si}
to give explicit formulas for elements of the quantum Mishchenko--Fomenko subalgebras
associated with the general linear and orthosymplectic Lie superalgebras (Propositions~\ref{prop:qmf}
and \ref{prop:qmfa}).

\section{Main theorem}
\label{sec:mth}

\subsection{Notation and definitions}
\label{subsec:def}

To introduce a suitable realization of the orthosymplectic Lie superalgebra
$\osp_{M|2n}$, we use
the
involution $i\mapsto i\pr=M+2n-i+1$ on
the set $\{1,2,\dots,M+2n\}$. Set
\ben
\bi=\begin{cases} 1\qquad\text{for}\quad i=1,\dots,n,n',\dots,1',\\
0\qquad\text{for}\quad i=n+1,\dots,(n+1)\pr
\end{cases}
\een
and
\ben
\ta_i=\begin{cases} \phantom{-}1\qquad\text{for}\quad i=1,\dots,M+n,\\
-1\qquad\text{for}\quad i=M+n+1,\dots,M+2n.
\end{cases}
\een
A standard basis of the general linear Lie superalgebra $\gl_{M|2n}$ is formed by elements $E_{ij}$
of the parity $\bi+\bj\mod 2$ for $1\leqslant i,j\leqslant M+2n$, with the commutation relations
\beql{glndr}
[E_{ij},E_{kl}]
=\de_{kj}\ts E_{i\tss l}-\de_{i\tss l}\ts E_{kj}(-1)^{(\bi+\bj)(\bk+\bl)}.
\eeq
We will regard the orthosymplectic Lie superalgebra $\osp_{M|2n}$
as the subalgebra
of $\gl_{M|2n}$ spanned by the elements\footnote{\label{foot:primes}We follow the settings of \cite{mr:rb};
an isomorphism with the presentation
of \cite{mn:ss} is given by $F_{ij}\mapsto F_{i'j'}$.}
\beql{fij}
F_{ij}=E_{ij}-E_{j'i'}(-1)^{\bi\tss\bj+\bi}\ts\ta_i\ta_j.
\eeq
The Lie superalgebra
$\osp_{M|2n}$ has Lie superbracket
\beql{comrel}
\big[F_{ij}, F_{kl}\big]= \de_{kj}F_{il}
- \de_{il}F_{kj}(-1)^{(\bi+\bj)(\bk+\bl)}
-  \de_{ki'} F_{j'l}(-1)^{\bi\tss\bj+\bi} \ta_i\ta_j
+\de_{j'l} F_{ki'}(-1)^{\bi+\bj+\bi\bk +\bj\bk}\ta_i\ta_j,
\eeq
and the symmetry relation
\beql{sym}
F_{ij}=-F_{j'i'}(-1)^{\bi\tss\bj+\bi}\ta_i\ta_j.
\eeq

Consider the $\ZZ_2$-graded vector space $\CC^{M|2n}$ over $\CC$ with the basis
$e_1,e_2,\dots,e_{1'}$, where the parity of the basis vector
$e_i$ is defined to be $\bi\mod 2$.
Accordingly, equip
the endomorphism algebra $\End\CC^{M|2n}$ with a $\ZZ_2$-gradation, such that
the parity of the matrix unit $e_{ij}$ is given by
$\bi+\bj\mod 2$.
We will consider even square matrices with entries in $\ZZ_2$-graded algebras; their
$(i,j)$ entries will have the parity $\bi+\bj\mod 2$.
Such a matrix superalgebra with matrix entries in
a superalgebra $\Ac$ will be identified with the tensor product superalgebra
$\End\CC^{M|2n}\ot\Ac$, so that a matrix $A=[A_{ij}]$ is regarded as the element
\ben
A=\sum_{i,j=1}^{1'}e_{ij}\ot A_{ij}(-1)^{\bi\tss\bj+\bj}\in \End\CC^{M|2n}\ot\Ac.
\een
We will use the involutive matrix {\em super-transposition} $t$ defined by
$(A^t)_{ij}=A_{j'i'}(-1)^{\bi\bj+\bj}\tss\ta_i\ta_j$.
We will also regard $t$ as the linear map
\beql{suptra}
t:\End\CC^{M|2n}\to \End\CC^{M|2n}, \qquad
e_{ij}\mapsto e_{j'i'}(-1)^{\bi\bj+\bi}\tss\ta_i\ta_j.
\eeq
In the case of multiple tensor products of the endomorphism algebras,
we will indicate by $t_l$ the map \eqref{suptra}
acting on the $l$-th copy of $\End\CC^{M|2n}$.

Introduce the permutation operator $P$ by
\ben
P=\sum_{i,j=1}^{1'} e_{ij}\ot e_{ji}(-1)^{\bj}\in \End\CC^{M|2n}\ot\End\CC^{M|2n}
\een
and set
\ben
Q=P^{\tss t_1}=P^{\tss t_2}=\sum_{i,j=1}^{1'} e_{ij}\ot e_{i'j'}(-1)^{\bi\bj}\ts\ta_i\ta_j
\in \End\CC^{M|2n}\ot\End\CC^{M|2n}.
\een

By introducing the matrix $\hF=[\hF_{ij}]$ with $\hF_{ij}=F_{ij}(-1)^{\bi}$, we can write the
commutation relations \eqref{comrel} in a matrix form as
\beql{comrellie}
[\hF_1,\hF_2]=(P-Q)\hF_2-\hF_2(P-Q),
\eeq
where both sides are elements of the tensor product superalgebra
\ben
\End\CC^{M|2n}\ot\End\CC^{M|2n}\ot \U(\osp_{M|2n})
\een
with
\ben
\hF_1=\sum_{i,j=1}^{1'} e_{ij}\ot 1\ot \hF_{ij}(-1)^{\bi\tss\bj+\bj}\fand
\hF_2=\sum_{i,j=1}^{1'} 1\ot e_{ij}\ot \hF_{ij}(-1)^{\bi\tss\bj+\bj},
\een
while $P$ and $Q$ are identified with $P\ot 1$ and $Q\ot 1$, respectively.
Similarly, the symmetry relation \eqref{sym} can be written as
\beql{symma}
\hF+\hF^{\ts t}=0.
\eeq

\subsection{Segal--Sugawara vectors}
\label{subsec:ssv}

As in the Introduction, we set $\g=\osp_{M|2n}$ and use the notation $X[r]=Xt^r$
for elements of $\wh\g$ with $X\in\g$. In the extended Lie superalgebra
$\wh\g\oplus \CC \tau$, we have $[\tau, X[r]]=-r X[r-1]$ for $X\in\g$, and $\tau$ commutes with $K$.
Denote by
$\U$ the universal enveloping algebra $\U(\wh\g\oplus \CC \tau)$,
and consider the tensor product superalgebra
\beql{tenprka}
\underbrace{\End\CC^{M|2n}\ot\dots\ot\End\CC^{M|2n}}_k{}\otimes \U.
\eeq
For $1 \leqslant a \leqslant k$ and $r \in \ZZ$, define
\beql{fra}
\hF[r]_a \coloneqq \sum_{i,j = 1}^{M+2n} 1 ^{\otimes (a-1)} \otimes e_{ij} \otimes 1^{\otimes (k-a)} \otimes \hF_{ij}[r] (-1)^{\bi \bj+\bj}.
\eeq
The symmetric group $\Sym_k$ acts on the superspace $(\CC^{M|2n})^{\ot k}$
by permutations of the
tensor factors. Such permutations
can be naturally identified with elements of the tensor product of the endomorphism
superalgebras in \eqref{tenprka}.
Denote by $H^{(k)}$ the element of the algebra
\eqref{tenprka} (with the identity component in $\U$) that is
the image of the symmetriser $h^{(k)}\in\CC\Sym_k$
defined by
\beql{ha}
h^{(k)}=\frac{1}{k!}\sum_{s\in\Sym_k} s
\eeq
under the action of $\Sym_k$.
Furthermore,
let $\la=(\la_1,\dots,\la_{\ell})$ be a partition of $k$ of length $\ell=\ell(\la)$,
so that $\la_1\geqslant\dots\geqslant\la_{\ell}>0$ and $\la_1+\dots+\la_{\ell}=k$.
We denote by $c^{}_{\la}$ the number
of permutations in the symmetric group $\Sym_k$ of cycle type $\la$.
Set
\beql{fla}
\hF[-\la]=\hF[-\la_1]_1 \dots \hF[-\la_{\ell}]_{\ell},
\eeq
and for positive integers $k$ introduce elements $\Phi_k\in \U(t^{-1}\g[t^{-1}])$ by
\beql{phimliesup}
\Phi^{}_{k}=\sum_{\la\ts\vdash k,\  \ell(\la)\ts \text{even}}\ts \Yc_{k,\ell}(M-2n-1)\ts
c^{}_{\la}\ts \str^{}_{1,\dots,\ell}\ts H^{(\ell)} \hF[-\la],
\eeq
noting that $\Phi_1=0$.
Here we use the polynomials $\Yc_{k,\ell}(T)$ in a variable $T$ defined by
\ben
\Yc_{k,\ell}(T)=\frac{\ell!}{k!}\ts\prod_{i=\ell}^{k-1}(T+i),
\een
while the supertrace
\ben
\str:\End \CC^{M|2n}\to\CC,\qquad e_{ij}\mapsto \de_{ij}(-1)^{\bi},
\een
is taken over the first $\ell$ copies of $\End \CC^{M|2n}$.
By the main theorem of \cite{mn:ss},
all elements $\Phi_k$ are {\em Segal--Sugawara vectors} associated with $\g$; that is,
they belong to the Feigin--Frenkel centre $\z(\wh\g)$.

\subsection{Harish-Chandra homomorphism}
\label{subsec:hch}

Consider the triangular decomposition
\beql{tridec}
\g=\n_-\oplus\h\oplus \n_+
\eeq
of the Lie superalgebra $\g$, where the subalgebras $\n_-$ and $\n_+$
are spanned by the elements $F_{ij}$ with $i>j$ and $i<j$, respectively, whereas
$\h$ is spanned by the elements $F_{ii}$ for $i=1,\dots,1'$.
The adjoint action of $\h$ on $t^{-1}\g[t^{-1}]$
extends to the universal enveloping algebra and
we have the homomorphism
for the $\h$-centralizer,
\beql{hchaff}
\f:\U\big(t^{-1}\g[t^{-1}]\big)^{\h}\to \U\big(t^{-1}\h[t^{-1}]\big),
\eeq
which is the projection to the first summand in the direct sum decomposition
\ben
\U\big(t^{-1}\g[t^{-1}]\big)^{\h}= \U\big(t^{-1}\h[t^{-1}]\big)\oplus
\Big(\U\big(t^{-1}\g[t^{-1}]\big)^{\h}\cap \U\big(t^{-1}\g[t^{-1}]\big)\ts t^{-1}\n_-[t^{-1}]\Big).
\een
The second summand is the kernel of the projection and coincides with
\ben
\U\big(t^{-1}\g[t^{-1}]\big)^{\h}\cap t^{-1}\n_+[t^{-1}]\ts \U\big(t^{-1}\g[t^{-1}]\big).
\een
We can regard $\f$ as an affine super-version of the {\em Harish-Chandra homomorphism};
cf. \cite[Sec.~2.2]{cw:dr} and \cite[Sec.~13.1]{m:so}.
Note that the Feigin--Frenkel centre $\z(\wh\g)$ is contained
in the centralizer $\U(t^{-1}\g[t^{-1}])^{\h}$.
To state our main result, set
\beql{muir}
\mu_i[-r]=F_{ii}[-r],\qquad i=1,\dots,1'\fand r=1,2,\dots,
\eeq
and assume that $M\geqslant 1$.
Suppose that $x_1,x_2,\dots,x_{1'}$ are variables (generally, noncommutative),
regarded as
elements of an associative algebra. For any nonnegative integer $k$,
introduce a noncommutative version of complete supersymmetric functions
by
\beql{nonchl}
h_k(x_1,x_2,\dots,x_{1'})=\sum_{a_1+\dots+a_{1'}=k}
x^{a_1}_{1}\dots x^{a_{1'}}_{1'},
\eeq
where $a_1,\dots,a_{1'}$ run over nonnegative integers and each of $a_i,a_{i'}$
takes only two values $0$ and $1$ for $i=1,\dots,n$. In the formulas below,
it is understood that $\tau\ts 1=0$.

\bth\label{thm:hch}
If $M=2m+1$, then the Harish-Chandra image $\f(\Phi_k)$ equals
\begin{multline}\non
h_k\big({-}\tau+\mu_{1}[-1],\dots,-\tau+\mu_{n}[-1],\tau+\mu_{n+1}[-1],\dots,\tau+\mu_{n+m}[-1],\\[0.3em]
\tau-\mu_{n+m}[-1],\dots,\tau-\mu_{n+1}[-1],
-\tau-\mu_{n}[-1],\dots,-\tau-\mu_{1}[-1]\big)\ts 1.
\end{multline}
\par
If $M=2m$, then the Harish-Chandra image $\f(\Phi_k)$ equals
the half-sum of two
expressions
\begin{multline}\non
h_k\big({-}\tau+\mu_{1}[-1],\dots,-\tau+\mu_{n}[-1],\tau+\mu_{n+1}[-1],\dots,\tau+\mu_{n+m-1}[-1],\\[0.3em]
\tau-\mu_{n+m}[-1],\dots,\tau-\mu_{n+1}[-1],
-\tau-\mu_{n}[-1],\dots,-\tau-\mu_{1}[-1]\big)\ts 1
\end{multline}
and
\begin{multline}\non
h_k\big({-}\tau+\mu_{1}[-1],\dots,-\tau+\mu_{n}[-1],\tau+\mu_{n+1}[-1],\dots,\tau+\mu_{n+m}[-1],\\[0.3em]
\tau-\mu_{n+m-1}[-1],\dots,\tau-\mu_{n+1}[-1],
-\tau-\mu_{n}[-1],\dots,-\tau-\mu_{1}[-1]\big)\ts 1.
\end{multline}
\eth

A proof of Theorem~\ref{thm:hch} will be given in Section~\ref{sec:pt}, following some applications
to Casimir elements, Sugawara operators, higher Gaudin Hamiltonians and shift-of-argument
subalgebras considered in Sections~\ref{sec:ce}--\ref{sec:gh}.

\subsection{Affine supersymmetric functions}
\label{subsec:asf}

With recent progress in proving the conjectures
concerning the Feigin--Frenkel centre $\z(\wh\gl_{m|n})$
in \cite{afn:ca}, including a connection
with the algebra of affine symmetric functions, it seems natural to expect
a similar connection in the orthosymplectic case.

To state the corresponding conjectures, follow \cite[Sec.~2.2]{cw:dr} to
recall the description of
the subalgebras of invariant polynomials $\La_{M|2n}(\mu)$ in the symmetric algebra $\Sr(\h)$
in the context of the Harish-Chandra isomorphism.
We will use the basis $\mu_1,\dots,\mu_{n+m}$ of the Cartan subalgebra $\h$ with $\mu_i=F_{ii}$.
If $M=2m+1$, then $\La_{2m+1|2n}(\mu)$
is the algebra of {\em supersymmetric polynomials} in the two sets
of variables $\mu_1^2,\dots,\mu_n^2$ and $\mu_{n+1}^2,\dots,\mu_{n+m}^2$.
Such a polynomial $P$ is symmetric in each set of variables and has the {\em cancellation
property}: the result of setting $\mu_1^2=\mu_{n+1}^2=z$ in $P$ is independent of $z$;
see e.g. \cite[Sec.~A.2]{cw:dr}.

The algebra $\La_{2m|2n}(\mu)$ is generated by all supersymmetric polynomials in the sets of
variables $\mu_1^2,\dots,\mu_n^2$ and $\mu_{n+1}^2,\dots,\mu_{n+m}^2$, and
the polynomial
\beql{suppf}
\mu_{n+1}\dots\mu_{n+m}\ts\prod_{i=1}^n\prod_{j=1}^m(\mu_i^2-\mu_{n+j}^2).
\eeq

To define affine supersymmetric functions associated with $\osp_{M|2n}$,
identify $\Sr(\h)$ with a subalgebra of the symmetric algebra $\Sr(t^{-1}\h[t^{-1}])$
by setting $\mu_i\mapsto\mu_i[-1]$.
The algebra $\Sr(t^{-1}\h[t^{-1}])$ is equipped with the derivation
$\tau$ acting on the generators by
\ben
\tau:\mu_i[-r]\mapsto r\tss\mu_i[-r-1],\qquad i=1,\dots,n+m,\qquad r=1,2,\dots.
\een

\bde\label{def:aff}
The {\em algebra $\La^{\aff}_{M|2n}(\mu)$} is defined as the subalgebra of
$\Sr(t^{-1}\h[t^{-1}])$ generated by all elements of the form $\tau^k P$,
where $k=0,1,\dots$ and $P\in \La_{M|2n}(\mu)$.
\qed
\ede

It is clear that the algebra $\La^{\aff}_{M|2n}(\mu)$ is also generated by
all elements of the form $\tau^k P$, where $P$ runs over a set of generators of $\La_{M|2n}(\mu)$.

Using the vector space isomorphism
$V_{\cri}(\g)\cong\U(t^{-1}\g\tss[t^{-1}])$,
equip $V_{\cri}(\g)$ with an ascending filtration induced by
the canonical filtration on the universal enveloping algebra.
The graded vector space $\gr V_{\cri}(\g)$ can therefore be identified
with the symmetric algebra $\Sr(t^{-1}\g[t^{-1}])$, while
the graded image of the Feigin--Frenkel centre $\z(\wh\g)$ is contained in the subalgebra
of $\g[t]$-invariants in $\Sr(t^{-1}\g[t^{-1}])$.
A super-commutative counterpart of \eqref{hchaff} is
the {\em Chevalley projection}
\beql{affche}
\bar\f:\Sr(t^{-1}\g\tss[t^{-1}])\to \Sr(t^{-1}\h\tss[t^{-1}])
\eeq
modulo the ideal $\Sr(t^{-1}\g\tss[t^{-1}])\big(t^{-1}\n_-[t^{-1}]\cup t^{-1}\n_+[t^{-1}]\big)$.
The formulas of Theorem~\ref{thm:hch} suggest the following
orthosymplectic version of
\cite[Conjecture~4.3]{mm:iv} concerning the Lie superalgebra $\gl_{m|n}$;
its proof for $\gl_{m|1}$ is given in \cite{afn:ca}.

\bcj\label{conj:chevalleyosp}
The restriction of the map \eqref{affche} to the subalgebra of $\g[t]$-invariants
yields an isomorphism of graded algebras
\ben
\Sr\big(t^{-1}\g[t^{-1}]\big)^{\g[t]}\cong \La^{\aff}_{M|2n}(\mu).
\een
\ecj

By the Harish-Chandra isomorphism for the Lie superalgebra $\osp_{2m|2n}$
(see \cite[Sec.~2.2]{cw:dr}), there exists a distinguished element
of the centre of $\U(\osp_{2m|2n})$
(the `noncommutative super-Pfaffian'),
corresponding to the polynomial
\eqref{suppf}, although its explicit form appears to be unknown.
Similarly, one can expect that there exists a distinguished Segal--Sugawara vector
$\Phi'\in\z(\wh\osp_{2m|2n})$, analogous to the Pfaffian-type vector in $\z(\wh\oa_{2m})$;
see \cite[Sec.~8.1]{m:so}. With this assumption, we can state the following conjecture
for the case of even $M$; cf. \cite[Conjecture~2.3]{mn:ss} for odd $M$.

\bcj\label{conj:even}
The elements $\Phi^{}_{2}, \Phi^{}_{4}, \Phi^{}_{6},\dots$ together with $\Phi'$
generate $\z(\wh\osp_{2m|2n})$ as a differential superalgebra.
\qed
\ecj

\section{Casimir elements}
\label{sec:ce}

For any nonzero $z\in\CC$ we have the evaluation homomorphism
\beql{evalrgm}
\ev_z:\U\big(t^{-1}\g[t^{-1}]\big)\to \U(\g),
\qquad F_{ij}[-r]\mapsto F_{ij}\tss z^{-r},\quad r>0.
\eeq
It is easy to see that the image of the Feigin--Frenkel centre $\z(\wh\g)$
under this homomorphism is contained in the centre $\Zr(\g)$ of
the universal enveloping algebra $\U(\g)$; cf. \cite[Prop.~6.5.2]{m:so}.

Following \cite{mn:ss}, note that the Segal--Sugawara vectors \eqref{phimliesup} can be defined
by an equivalent formula involving
the symmetrizer $S^{(k)}$ in the Brauer algebra $\Bc_k(M-2n)$
as recalled in \eqref{symo} below.
If $M$ is odd, then
the vector $\Phi_k$ coincides with the constant term of the polynomial in $\tau$
given by
\beql{sspol}
\ga_k(M-2n)\ts\str^{}_{1,\dots,k}\tss S^{(k)} \big(\tau+\hF[-1]_{1}\big)\dots
\big(\tau+\hF[-1]_{k}\big),
\eeq
where
\beql{ga}
\ga_k(\om)=\frac{\om+k-2}{\om+2\tss k-2}.
\eeq
According to \cite[Remark~4.9(i)]{mn:ss}, the polynomial \eqref{sspol} expands as
\beql{poltaexpa}
\sum_{r=0}^k \ts\binom{M-2n+k-2}{k-r}\ts \Phi_r\ts \tau^{k-r},
\eeq
with $\Phi_0=1$.
The same is true for even values $M=2m$ only if $m>n$. However, the arguments
of Sec.~\ref{subsec:ea} imply
that the polynomial \eqref{sspol}, regarded as a function of $m$ with the parameters $n$ and $k$ fixed,
has a removable singularity at the points where $m\leqslant n$.
The values of the coefficients of the polynomial \eqref{sspol}
at these points can therefore be considered to be equal to
the respective coefficients in \eqref{poltaexpa}.
In what follows, we will use expression \eqref{sspol} for all values
of the parameters, assuming this interpretation.

The image of the polynomial \eqref{sspol} under the evaluation homomorphism \eqref{evalrgm}
is given by
\beql{sspolev}
\ga_k(M-2n)\ts\str^{}_{1,\dots,k}\tss S^{(k)} \big({-}\di_z+\hF_{1}\tss z^{-1}\big)\dots
\big({-}\di_z+\hF_{k}\tss z^{-1}\big),
\eeq
where $-\di_z$ is understood as the image of $\tau$.
Multiply this expression by $z^k$ from the left
and use the relations
\ben
z^l\big({-}\di_z+\hF_{i}\tss z^{-1}\big)=\big({-}z\di_z+l-1+\hF_{i}\big)\tss z^{l-1}
\een
to get a polynomial in $u$ with coefficients in the center $\Zr(\g)$,
\beql{sspolevu}
\ga_k(M-2n)\ts\str^{}_{1,\dots,k}\tss S^{(k)} \big(u+k-1+\hF_{1}\big)\dots
\big(u+\hF_{k}\big),
\eeq
where we set $u=-z\di_z$. By using \eqref{poltaexpa}, we derive
an equivalent formula for this polynomial,
\beql{eqpolu}
\sum_{r=0}^{k}\binom{k}{r}\ts u^{k-r}\ts
\sum_{\la\ts\vdash r,\  \ell(\la)\ts \text{even}}\ts \Yc_{k,\ell}(M-2n-1)\ts
c^{}_{\la}\ts \str^{}_{1,\dots,\ell}\ts H^{(\ell)} \hF_1\dots \hF_{\ell},
\eeq
where the second sum is understood as equal to $\Yc_{k,0}(M-2n-1)$ for $r=0$.
Our goal is to use this polynomial to produce
generators of $\Zr(\g)$ and describe them by calculating their images under
the Harish-Chandra isomorphism.

To recall the construction of the isomorphism,
use the triangular decomposition
\eqref{tridec}
and note that the adjoint action of $\h$ on $\g$
extends to the universal enveloping algebra. Hence
we have the homomorphism
for the $\h$-centralizer,
\beql{hchfd}
\chi:\U(\g)^{\h}\to \U(\h)
\eeq
which is the projection to the first summand in the direct sum decomposition
\ben
\U(\g)^{\h}= \U(\h)\oplus
\big(\U(\g)^{\h}\cap \U(\g)\ts \n_+\big).
\een
The {\em Harish-Chandra isomorphism}
\beql{chi}
\chi:\Zr(\g)\to \La_{M|2n}(\ell)
\eeq
is obtained by
the restriction of $\chi$ to the centre $\Zr(\g)$; see e.g. \cite[Sec.~2.2]{cw:dr}.
To describe the image $\La_{M|2n}(\ell)\subset \U(\h)$,
introduce the notation
$\la_i=F_{ii}$ for $i=1,\dots,1'$, so that $\la_{i'}=-\la_i$.
Set
\beql{lairfd}
\ell_i=\la_{i}+\rho_i,\qquad i=1,\dots,n+m,
\eeq
where $\rho=(\rho_1,\dots,\rho_{n+m})$ is the {\em Weyl vector}
which is the difference between the half-sums of positive even and odd roots,
$\rho=\rho^{}_{\bar 0}-\rho^{}_{\bar 1}$. Explicitly, the components of $\rho$
are given by
\ben
\rho_i=-\frac{M}{2}+n-i+1,\quad i=1,\dots,n,\Fand \rho_{n+j}=\frac{M}{2}-j,\quad j=1,\dots,m.
\een
As in Sec.~\ref{subsec:asf},
if $M=2m+1$, then $\La_{2m+1|2n}(\ell)$
is the algebra of supersymmetric polynomials in two sets
of variables $\ell_1^2,\dots,\ell_n^2$ and $\ell_{n+1}^2,\dots,\ell_{n+m}^2$.
The algebra $\La_{2m|2n}(\ell)$ is generated by
all supersymmetric polynomials in the sets of
variables $\ell_1^2,\dots,\ell_n^2$ and $\ell_{n+1}^2,\dots,\ell_{n+m}^2$, and
the polynomial
\beql{suppfell}
\pi=\ell_{n+1}\dots\ell_{n+m}\ts\prod_{i=1}^n\prod_{j=1}^m(\ell_i^2-\ell_{n+j}^2).
\eeq
Since $\pi^2$ is a supersymmetric polynomial, the elements of $\La_{2m|2n}(\ell)$
are polynomials of the form $P_1+\pi\tss P_2$, where $P_1$ and $P_2$
are supersymmetric.

By an equivalent viewpoint, the isomorphism \eqref{chi} can be defined
using {\em highest weight $\g$-modules} $L(\la)$. Such a module
with highest weight $\la=(\la_1,\dots,\la_{n+m})$ is generated
by
a nonzero vector $v\in L(\la)$
such that
\begin{alignat}{2}
F_{ij}\ts v&=0 \qquad &&\text{for} \quad
1\leqslant i<j\leqslant 1', \qquad \text{and}
\non
\\
F_{ii}\ts v&=\la_i\ts v \qquad &&\text{for} \quad 1\leqslant i\leqslant n+m,
\non
\end{alignat}
where the $\la_i$ are now understood as complex numbers.
Any element $z\in\Z(\g)$ acts in $L(\la)$
by multiplying each vector by a scalar $\chi(z)$.
When regarded as a function of the highest weight, $\chi(z)$
is a polynomial in $\la_1,\dots,\la_{n+m}$ that
belongs to the algebra $\La_{M|2n}(\ell)$.

To apply Theorem~\ref{thm:hch}, note that the roles
of the subalgebras $\n_+$ and $\n_-$ in the definitions
of the homomorphisms \eqref{hchaff} and \eqref{hchfd} are swapped.
To be able to use the homomorphism \eqref{evalrgm},
observe that the linear map taking $F_{ij}$ to $F_{i'j'}\ta_i\ta_j=-F_{ji}(-1)^{\bi\bj+\bj}$
defines an automorphism of $\g$.
The composition of the action of $\g$ on the highest weight module $L(\la)$
with this automorphism
yields a {\em lowest weight module} with the lowest weight $-\la$, because under the new action,
the vector $v$ is annihilated by the subalgebra $\n_-$.
Thus, we arrive at the following corollary of Theorem~\ref{thm:hch}.

\bco\label{cor:hchfd}
If $M=2m+1$, then the Harish-Chandra image of the polynomial \eqref{sspolevu} equals
\begin{multline}\non
(-1)^k\ \sum_{i_1<\dots<i_s\leqslant n<i_{s+1}\leqslant\dots\leqslant i_p<n'\leqslant i_{p+1}<\dots<i_k}
(u+k-1+\la_{i_1})\dots(u+k-s+\la_{i_s})\\[0.4em]
{}\times (-u-k+s+1+\la_{i_{s+1}})\dots
(-u-k+p+\la_{i_{p}})(u+k-p-1+\la_{i_{p+1}})\dots (u+\la_{i_k}),
\end{multline}
assuming varying values of the parameters
$s$ and $p$ in the summation, and $n+m+1$ does not
occur among the summation indices.
\par
If $M=2m$, then the Harish-Chandra image of the polynomial \eqref{sspolevu} equals
the half-sum of two above expressions, where $n+m$ does not
occur among the summation indices in the first expression, while $(n+m)'$
does not occur in the second.
\eco

\bpf
As with the calculation of the polynomial \eqref{sspolevu},
multiply the
Harish-Chandra images in Theorem~\ref{thm:hch} by $z^k$ from the left
and replace $\mu_i[-1]$ with $-\la_i\tss z^{-1}$.
\epf

Due to the isomorphism \eqref{chi}, all coefficients of the polynomials in Corollary~\ref{cor:hchfd}
belong to the subalgebra $\La_{M|2n}(\ell)$. We will choose suitable values of $u$ and $k$
to produce generators of the centre $\Zr(\g)$ and calculate their images in $\La_{M|2n}(\ell)$.
Namely, we only take even values of $k$ and so replace $k$ with $2k$, and set
$u=-k$ (taking $u=-k+1$ would also give the same elements).
For all $k=1,2,\dots$ introduce the elements
\beql{ck}
C_k=\ga^{}_{2k}(M-2n)\ts\str^{}_{1,\dots,2k}\tss S^{(2k)} \big(\hF_{1}+k-1\big)\dots
\big(\hF_{2k}-k\big).
\eeq
More explicitly, by using \eqref{eqpolu},
we get an equivalent expression
\ben
C_k=\sum_{r=0}^{2k}\binom{2k}{r}(-k)^{2\tss k-r}\ts
\sum_{\la\ts\vdash r,\  \ell(\la)\ts \text{even}}\ts \Yc_{2k,\ell}(M-2n-1)\ts
c^{}_{\la}\ts \str^{}_{1,\dots,\ell}\ts H^{(\ell)} \hF_1\dots \hF_{\ell}.
\een
Also, denote by $C^{\tss\prime}$ the Casimir element in $\Zr(\osp_{2m|2n})$ whose Harish-Chandra image
coincides with the polynomial
\eqref{suppfell}.

To state the main result of this section, assume that $M\geqslant 1$ and recall the {\em complete factorial supersymmetric polynomials}
$h_k(x/y\, |\,a)$ as introduced in \cite{m:fs}. We consider two sets of variables $x=(x_1,\dots,x_m)$
and $y=(y_1,\dots,y_n)$ along with a parameter sequence $a=(a_i)$, $i\in\ZZ$. Then
the polynomials are defined by
\begin{multline}
h_k(x/y\, |\,a)={}\\[0.4em]
 \sum_{p+q=k}\sum_{\begin{smallmatrix}i_1\leqslant \dots \leqslant i_p \\ j_1>\cdots > j_q\end{smallmatrix}} (-1)^q\tss(y_{j_1}- a_{j_1}) \cdots (y_{j_q} - a_{j_q + q-1}) (x_{i_1} - a_{i_1+q}) \cdots (x_{i_p} - a_{i_p+k-1}).
\label{hkxya}
\end{multline}
They are supersymmetric in $x$ and $y$ in the sense of the definition
recalled in Sec.~\ref{subsec:asf}; note that our variables $y_i$ are related to
those in \cite{m:fs} by the change of signs $y_i\to -y_i$ for $i=1,\dots,n$.

\bth\label{thm:hchfd}
\ \ (i)\quad
The elements $C_{1}, C_{2},\dots$ generate the centre $\Zr(\osp_{2m+1|2n})$.
\par
(ii)\quad The elements $C_{1}, C_{2},\dots$ together with $C^{\tss\prime}$ generate
the centre $\Zr(\osp_{2m|2n})$.
\par
(iii) \quad In both cases, the Harish-Chandra images of the generators $C_k$ are given by
\beql{chick}
\chi(C_k)=h_k(x/y\, |\,a),
\eeq
where the variables are specialized by $x=(\ell^2_{n+1},\dots,\ell^2_{n+m})$ and
$y=(\ell^2_{1},\dots,\ell^2_{n})$, while the components of the
sequence $a$ are $a_i=\big(\frac{M}2-i+k-1\big)^2$.
\eth

\bpf
Since the polynomials $h_k(x/y\, |\,a)$ generate the algebra of supersymmetric polynomials,
it is enough to prove part (iii). Suppose first that $M=2m+1$.
The Harish-Chandra image
$\chi(C_k)$ is given by Corollary~\ref{cor:hchfd} and it has the form
\begin{multline}\label{polla}
\sum_{i_1<\dots<i_s\leqslant n<i_{s+1}\leqslant\dots\leqslant i_p<n'\leqslant i_{p+1}<\dots<i_{2k}}
(\la_{i_1}+k-1)\dots(\la_{i_s}+k-s)\\[0.4em]
{}\times (\la_{i_{s+1}}-k+s+1)\dots
(\la_{i_{p}}-k+p)(\la_{i_{p+1}}+k-p-1)\dots (\la_{i_{2k}}-k).
\end{multline}
To show that this coincides with $h_k(x/y\, |\,a)$, we
use induction on $n$ by taking $n=0$ as the induction base.
In this particular case,
Theorem~\ref{thm:hchfd} reduces to the corresponding theorem
for the orthogonal Lie algebras $\oa_M$ in \cite[Prop.~5.4.1 and Sec.~13.4]{m:so}
and holds trivially for $m=0$.
The Harish-Chandra image in \eqref{chick} becomes
the {\em complete factorial symmetric polynomial} $h_k(x\, |\,a)$.
It is easy to verify by using the symmetry of the polynomials in the variables $x_1,\dots,x_m$,
that $h_k(x\, |\,a)$ coincides with $h_k(x\, |\,a')$, where the sequence $a'=(a'_i)$
is defined by $a'_i=(i-\frac12)^2$, as used therein.

Now suppose that $n\geqslant 1$ and observe that
the polynomial \eqref{polla}, as an element of $\La_{M|2n}(\ell)$,
has degree $1$ with respect to the variable
$\ell_1^2$.
Under the evaluation $\la_1=-k+1$ we have $\ell_1^2=a_n$. On the other hand,
putting $\la_1=-k+1$
in \eqref{polla} yields a polynomial of the same form in $\la_2,\dots,\la_{n+m}$
with $n$ replaced by $n-1$.
Indeed, the terms with $i_{1}=1$ vanish, whereas the sum of the remaining terms with $i_{2k}=1'$
equals the negative of
\begin{multline}\non
\sum_{2\leqslant i_1<\dots<i_s\leqslant n<i_{s+1}\leqslant\dots\leqslant i_p<n'
\leqslant i_{p+1}<\dots<i_{2k-1}\leqslant 2'}
(\la_{i_1}+k-1)\dots(\la_{i_s}+k-s)\\[0.4em]
{}\times (\la_{i_{s+1}}-k+s+1)\dots
(\la_{i_{p}}-k+p)(\la_{i_{p+1}}+k-p-1)\dots (\la_{i_{2k-1}}-k+1).
\end{multline}
However, this sum is zero, which is seen by using the involution on the multiset of indices
taking
$(i_1,\dots,i_{2k-1})$ to $(i'_{2k-1},\dots,i'_{1})$, and recalling that $\la_{i'}=-\la_i$.
By the induction hypothesis, the result of the evaluation
equals $h_k(x/y'\, |\,a)$, where $y'=(\ell_2^2,\dots,\ell_n^2)$.
Furthermore, the induction hypothesis also implies that
the coefficient of $\la^2_1$ in \eqref{polla} coincides with $-h_{k-1}(x/y'\, |\,b)$,
where the components of the sequence $b$ are given by $b_i=a_{i+1}$.
Hence, the polynomial \eqref{polla} can be written as
\beql{hkexp}
h_k(x/y'\, |\,a)-(\ell_1^2-a_n)\tss h_{k-1}(x/y'\, |\,b),
\eeq
and we will show that this
coincides with $h_k(x/y\, |\,a)$. Indeed, write
$h_k(x/y\, |\,a)$ in an equivalent form by using the symmetry on the $y$ variables,
replacing $y_i\mapsto y_{n-i+1}$ to get
\begin{multline}
h_k(x/y\, |\,a)={}\\[0.4em]
 \sum_{p+q=k}\sum_{\begin{smallmatrix}i_1\leqslant \dots \leqslant i_p \\ r_1<\cdots < r_q\end{smallmatrix}} (-1)^q\tss(y_{r_1}- a_{n-r_1+1}) \cdots (y_{r_q} - a_{n-r_q + q}) (x_{i_1} - a_{i_1+q}) \cdots (x_{i_p} - a_{i_p+k-1}).
\non
\end{multline}
Split the sum by separating the terms with $r_1=1$ and replace $y_1$ by $\ell_1^2$ to see that
this coincides with \eqref{hkexp}, which completes the argument in the case under consideration.

In the case $M=2m$,
the Harish-Chandra image
$\chi(C_k)$ provided by Corollary~\ref{cor:hchfd} equals
the half-sum of two expressions given by \eqref{polla} with the respective
restrictions of the summation indices. This implies that the image
is stable under the transformation $\la_{n+m}\mapsto \la_{(n+m)'}$
which can be written as $\ell_{n+m}\mapsto -\ell_{n+m}$.
The polynomial $\pi$ defined in \eqref{suppfell} changes sign under this transformation
and therefore it cannot occur in the expression of
$\chi(C_k)$ as an element of $\La_{2m|2n}(\ell)$.
In other words, $\chi(C_k)$ is a supersymmetric polynomial
in the sets of variables $\ell_1^2,\dots,\ell_n^2$ and $\ell_{n+1}^2,\dots,\ell_{n+m}^2$.
The proof is completed by the same induction argument over $n$
as for the odd $M$,
where the base
case $n=0$ corresponds to the results
for the orthogonal Lie algebras in \cite[Prop.~5.4.1 and Sec.~13.4]{m:so}.
\epf

Note that in the special case $M=1$,
the generators $C_1,\dots,C_n$ of $\Zr(\osp_{1|2n})$
are algebraically independent since $\La_{1|2n}(\ell)$
is the algebra of symmetric polynomials in $\ell^2_{1},\dots,\ell^2_{n}$.

The case $M=0$ corresponds to the symplectic Lie algebras and it is not covered
by Theorem~\ref{thm:hchfd}. The Harish-Chandra images of the elements $C_k$ in \eqref{ck}
have a different form
as given in \cite[Prop.~5.5.4]{m:so}; a new proof is implied
by the results of Sec.~\ref{subsec:np} below.

Note also that
the Harish-Chandra image of the polynomial \eqref{sspolevu}
(for any values of $M$ and $n$)
can be calculated in terms of the
variables $\ell_i^2$ by using arguments similar to those in \cite[Sec.~5.4]{m:so}.

\section{Sugawara operators}
\label{sec:so}

Here we will use Theorem~\ref{thm:hch} to calculate the
Harish-Chandra images of the
{\em Sugawara operators} for the Lie superalgebra $\g=\osp_{M|2n}$. They are
elements of the centre $\Z(\wh\g)$ of the
completed universal enveloping algebra $\wt\U_{\cri}(\wh\g)$
at the critical level. We will follow the book \cite[Sections~3.2 and 4.3]{f:lc}
to recall the key definitions and connection between the Feigin--Frenkel centre $\z(\wh\g)$
as defined in \eqref{centsup}, and the algebra $\Z(\wh\g)$.

Consider the quotient $\U_{\cri}(\wh\g)$ of the
universal enveloping algebra $\U(\wh\g)$ by the
ideal generated by $K+h^\vee$,
where $h^\vee=M-2n-2$.
Introduce a linear topology on $\U_{\cri}(\wh\g)$ by using the
neighborhood basis for $0$ formed by the left ideals
$I_p$ of $\U_{\cri}(\wh\g)$ generated
by $t^p\g[t]$ for all $p\geqslant 0$.
The completed universal enveloping algebra $\wt\U_{\cri}(\wh\g)$
is the
completion of $\U_{\cri}(\wh\g)$ with respect to this topology.
Equivalently, $\wt\U_{\cri}(\wh\g)$ can be defined as the inverse limit
\beql{complua}
\wt\U_{\cri}(\wh\g)=\lim_{\longleftarrow}\U_{\cri}(\wh\g)/I_p.
\eeq

For any $i,j\in\{1,\dots,1'\}$ introduce the Laurent series $F_{ij}(u)$
with coefficients in $\U_{\cri}(\wh\g)$ by
\ben
F_{ij}(u)=\sum_{r\in\ZZ} F_{ij}[r]\tss u^{-r-1}.
\een
Combine these series into the matrix $\hF(u)=[\hF_{ij}(u)]$ with $\hF_{ij}(u)=F_{ij}(u)(-1)^{\bi}$.
The {\em state-field correspondence map}
\ben
Y:V_{\cri}(\g)\to \wt\U_{\cri}(\wh\g)
\een
is a linear map defined by
\ben
Y:F_{i_1j_1}[-r_1-1]\dots F_{i_kj_k}[-r_k-1]
\mapsto
\frac{1}{r_1!\ts\dots r_k!}: F^{(r_1)}_{i_1j_1}(u)\dots F^{(r_k)}_{i_kj_k}(u):
\een
for any nonnegative integers $r_1,\dots,r_k$, where $F^{(r)}_{ij}(u)=\di_u^{\ts r}F_{ij}(u)$
and we use the standard normal ordering notation.
Namely, the normally ordered product of homogeneous fields
\ben
a(u)=\sum_{r\in\ZZ}a_{(r)}u^{-r-1}\Fand
b(w)=\sum_{r\in\ZZ}b_{(r)}w^{-r-1}
\een
of the respective parities $\ba$ and $\bb$
is the formal power series
\beql{normor}
:a(u)\tss b(w)\ts{:}
= a(u)_+\tss b(w)+(-1)^{\ba\tss\bb}\ts b(w)\tss a(u)_-,
\eeq
where
\ben
a(u)_+=\sum_{r<0}a_{(r)}u^{-r-1}\Fand
a(u)_-=\sum_{r\geqslant 0}a_{(r)}u^{-r-1}.
\een
This definition extends to an arbitrary number of fields
with the convention that the normal ordering is read
from right to left.
The image of the restriction of the map $Y$ to the subspace $\z(\wh\g)\subset V_{\cri}(\g)$
is contained in the centre $\Z(\wh\g)$ of $\wt\U_{\cri}(\wh\g)$.

We extend the matrix notation \eqref{fra} to the superalgebras of the form
\eqref{tenprka}, where $\U$ is replaced by the completed universal enveloping algebra.
Denote by $\hF(u,\la)$ the image of the element $\hF[-\la]$ defined in
\eqref{fla} under the extension of the map $Y$ to the tensor superalgebra,
which acts as the identity map
on the endomorphism superalgebra factors.

Introduce the Laurent series
\beql{phiksuga}
\Phi^{}_{k}(u)=\sum_{\la\ts\vdash k,\  \ell(\la)\ts \text{even}}\ts \Yc_{k,\ell}(M-2n-1)\ts
c^{}_{\la}\ts \str^{}_{1,\dots,\ell}\ts H^{(\ell)} \hF(u,\la).
\eeq
The following proposition is immediate from \cite[Theorem~2.1]{mn:ss}.

\bpr\label{prop:suga}
All coefficients of the Laurent series
$\Phi^{}_{k}(u)$ for positive integers $k$
are Sugawara operators for $\osp_{M|2n}$.
\qed
\epr

Note that an equivalent formula for the series \eqref{phiksuga} can be written
with the use of the symmetrizer $S^{(k)}$ in the Brauer algebra $\Bc_k(M-2n)$
in the same way as
for the Segal--Sugawara vectors in \eqref{sspol} above; see the comments
following \eqref{poltaexpa}.
We have
\beql{phiksugasym}
\Phi^{}_{k}(u)=\ts :\tss\ga_k(M-2n)\ts\str^{}_{1,\dots,k}\tss S^{(k)} \big(\di_u+\hF(u)_{1}\big)\dots
\big(\di_u+\hF(u)_{k}\big)\tss 1 :\ts,
\eeq
assuming $\di_u\tss 1=0$, where $\ga_k(M-2n)$ is defined in \eqref{ga}.

Now recall the affine version of the Harish-Chandra homomorphism
involving the centre $\Z(\wh\g)$; cf. \cite[Sec.~13.3]{m:so}.
Keeping the triangular decomposition
\eqref{tridec} of $\g$, introduce any total ordering $\prec$ on the basis elements of the
affine Kac--Moody superalgebra $\wh\g$ defined in
\eqref{km} to satisfy the following conditions.
Each basis element of $t^{-1}\g[t^{-1}]$ should precede each basis element
of $\g[t]$, and
the ordering on the corresponding basis elements of $\wh\g$ should be consistent with
the conditions
\ben
\n_-[t]\prec \h[t]\prec \n_+[t]\Fand
t^{-1}\tss\n_+[t^{-1}]\prec t^{-1}\tss\h[t^{-1}]\prec t^{-1}\tss\n_-[t^{-1}],
\een
indicating the ordering between the basis elements belonging to the subspaces of $\wh\g$.

By the Poincar\'e--Birkhoff--Witt theorem,
any element $x\in\U(\wh\g)$
can be written as a unique linear combination of
ordered monomials
in the basis elements of $\wh\g$.
Set
\ben
\wh\h=\h\tss[t,t^{-1}]\oplus\CC K
\een
and denote by $x_0\in \U(\wh\h)$ the component
of the linear combination representing the element $x$,
where each
monomial does not
contain any basis elements $X[r]$ with $X\in\n_-\oplus\n_+$.
The linear map $\theta:x\mapsto x_0$ defines the projection
$\theta:\U(\wh\g)\to \U(\wh\h)$.
Extend $\theta$
by continuity to get the projection
\ben
\theta:\wt\U_{\cri}(\wh\g)\to \wt\U_{\cri}(\wh\h),
\een
where $\wt\U_{\cri}(\wh\h)$ denotes the completion
of $\U_{\cri}(\wh\h)$ at the critical level defined as in \eqref{complua}.

Note that, unlike the Cartan subalgebra $\h$, the Lie algebra $\wh\h$ is not abelian.
Therefore, we need one more step before restricting the projection $\theta$ to the centre $\Z(\wh\g)$.
We identify $\U_{\cri}(\wh\h)$ with the symmetric algebra
$\Pi:=\Sr(\h\tss[t,t^{-1}])$ via
the natural isomorphism of vector spaces by using the basis of
ordered monomials. It extends to an isomorphism of the respective completed
vector spaces
$
\eta:\wt\U_{\cri}(\wh\h)\to \wt\Pi$.
We thus get a linear map
$
\f:\wt\U_{\cri}(\wh\g)\to \wt\Pi
$
defined as the composition $\f=\eta\circ\theta$. By the same argument as in \cite[Prop.~13.3.1]{m:so},
the restriction of $\f$ to the centre $\Z(\wh\g)$ yields a homomorphism
\beql{hchcompl}
\f:\Z(\wh\g)\to \wt\Pi
\eeq
which can be regarded as an affine version of the
Harish-Chandra homomorphism.

The following corollary is immediate from Theorem~\ref{thm:hch}. Set
\ben
\mu_i[r]=F_{ii}[r]\Fand \mu_i(u)=\sum_{r\in\ZZ} \mu_i[r]\tss u^{-r-1}
\een
and use the notation \eqref{nonchl}.

\bco\label{cor:hchss}
If $M=2m+1$, then the Harish-Chandra image of the Laurent series $\Phi_k(u)$ equals
\begin{multline}\non
h_k\big({-}\di_u+\mu_{1}(u),\dots,-\di_u+\mu_{n}(u),\di_u+\mu_{n+1}(u),\dots,\di_u+\mu_{n+m}(u),\\[0.3em]
\di_u-\mu_{n+m}(u),\dots,\di_u-\mu_{n+1}(u),
-\di_u-\mu_{n}(u),\dots,-\di_u-\mu_{1}(u)\big).
\end{multline}
\par
If $M=2m$, then the Harish-Chandra image of the Laurent series $\Phi_k(u)$
equals
the half-sum of two
noncommutative complete symmetric functions
\begin{multline}\non
h_k\big({-}\di_u+\mu_{1}(u),\dots,-\di_u+\mu_{n}(u),\di_u+\mu_{n+1}(u),\dots,\di_u+\mu_{n+m-1}(u),\\[0.3em]
\di_u-\mu_{n+m}(u),\dots,\di_u-\mu_{n+1}(u),
-\di_u-\mu_{n}(u),\dots,-\di_u-\mu_{1}(u)\big)
\end{multline}
and
\begin{multline}\non
h_k\big({-}\di_u+\mu_{1}(u),\dots,-\di_u+\mu_{n}(u),\di_u+\mu_{n+1}(u),\dots,\di_u+\mu_{n+m}(u),\\[0.3em]
\di_u-\mu_{n+m-1}(u),\dots,\di_u-\mu_{n+1}(u),
-\di_u-\mu_{n}(u),\dots,-\di_u-\mu_{1}(u)\big).
\end{multline}
\eco

\section{Gaudin Hamiltonians and shift-of-argument subalgebras}
\label{sec:gh}

By a remarkable observation of Feigin, Frenkel and Reshetikhin~\cite{ffr:gm},
the centre of the affine vertex algebra at the critical level
is closely related to the Hamiltonians of the Gaudin model
describing quantum spin chains. Here we will apply this observation
to the orthosymplectic Lie superalgebras $\g=\osp_{M|2n}$ to produce
higher Hamiltonians associated with the Segal--Sugawara vectors $\Phi_k$.
Furthermore, we will follow Rybnikov~\cite{r:si} to produce a family
of commutative subalgebras of $\U(\g)$ parameterized by
functionals $\mu\in\g^*_{\bar 0}$. These are orthosymplectic versions $\Ac_{\mu}$ of the
{\em quantum shift-of-argument} or {\em Mishchenko--Fomenko subalgebras}
whose classical counterparts
originated in \cite{mf:ee}.

\subsection{Higher Hamiltonians}
\label{subsec:hh}

By the vacuum axiom of vertex algebras, the application of
the fields $\Phi_k(u)$
introduced in \eqref{phiksugasym}
to the vacuum vector yields formal power series in $u$
with coefficients in the Feigin--Frenkel centre $\z(\wh\g)$.
In particular, these coefficients
generate a commutative subalgebra
of $\U(t^{-1}\g[t^{-1}])$.
Explicit formulas for these coefficients are obtained
by replacing $\hF(u)$ with the matrix $\hF(u)_+=[\hF_{ij}(u)_+]$, where
\beql{fuplus}
\hF_{ij}(u)_+=\sum_{r=1}^{\infty} \hF_{ij}[-r]\tss u^{r-1},
\eeq
to get the series
\beql{gradepbd}
\Phi_k(u)_+=\ga_k(M-2n)\ts\str^{}_{1,\dots,k}\tss S^{(k)} \big(\di_u+\hF(u)_{+\ts 1}\big)\dots
\big(\di_u+\hF(u)_{+\ts k}\big)\tss 1.
\eeq
Given
a nonzero $z\in\CC$, consider the evaluation homomorphism \eqref{evalrgm}.
The image of the series \eqref{fuplus} is given by
\ben
\hF_{ij}(u)_+\mapsto \frac{\hF_{ij}}{z-u}.
\een

Using the coassociativity of the standard coproduct on $\U(t^{-1}\g[t^{-1}])$
defined by
\ben
\Delta: Y\mapsto Y\ot 1+1\ot Y,\qquad Y\in t^{-1}\g[t^{-1}],
\een
for any $\ell\geqslant 1$ we get the homomorphism
\beql{comult}
\U\big(t^{-1}\g[t^{-1}]\big)\to \U\big(t^{-1}\g[t^{-1}]\big)^{\ot\tss\ell}
\eeq
as an iterated coproduct map. Now fix distinct nonzero complex numbers
$z_1,\dots,z_\ell$ and apply homomorphisms of the form \eqref{evalrgm} to the tensor factors in
\eqref{comult}; we get another homomorphism
\beql{psiu}
\psi:\U\big(t^{-1}\g[t^{-1}]\big)\to \U(\g)^{\ot\tss\ell},
\eeq
so that
\ben
\psi:\hF_{ij}(u)_+\mapsto \sum_{a=1}^\ell \frac{(\hF_{ij})_a}{z_a-u}\in \U(\g)^{\ot\tss\ell},
\een
where $X_a=1^{\ot (a-1)}\ot X\ot 1^{\ot (\ell-a)}$ for $X\in \U(\g)$.
The image $\hF_{ij}(u)^{\psi}_+:=\psi(\hF_{ij}(u)_+)$ can be regarded as an operator
on the tensor product
\beql{gaudmo}
M_1\ot\dots\ot M_{\ell}
\eeq
of arbitrary $\g$-modules $M_1,\dots,M_{\ell}$, where $(\hF_{ij})_a$ is understood as an operator on $M_a$.
We thus arrive at the following conclusion.

\bpr\label{prop:gh}
The coefficients of all series
\ben
\ga_k(M-2n)\ts\str^{}_{1,\dots,k}\tss S^{(k)} \big(\di_u+\hF(u)^{\psi}_{+\ts 1}\big)\dots
\big(\di_u+\hF(u)^{\psi}_{+\ts k}\big)\tss 1
\een
with positive integer values of $k$ form a commutative family of operators
on the module \eqref{gaudmo}.
\qed
\epr

In the same way as for the simple Lie algebra case considered in \cite{ffr:gm},
the family of commuting operators provided by Proposition~\ref{prop:gh} contains
the quadratic Gaudin Hamiltonian. Namely,
note that $2\tss \Phi_2(u)_+=\str\ts \hF(u)_+^2$,
and the image of this series under the homomorphism \eqref{psiu} equals
\beql{hamt}
\sum_{a=1}^{\ell}\frac{\Hc^{(a)}}{z_a-u}+\sum_{a=1}^{\ell}\frac{C_a}{(z_a-u)^2},
\eeq
where
\ben
\Hc^{(a)}=\sum_{b\ne a}\frac{1}{z_b-z_a}\sum_{i,j=1}^{1'} (F_{ij})_a(F_{ji})_b(-1)^{\bj}
\een
and
\ben
C=\sum_{i,j=1}^{1'} F_{ij}F_{ji}(-1)^{\bj}.
\een
These formulas for the Hamiltonian \eqref{hamt} already appeared in \cite{km:co}
in the simplest case of $\osp_{1|2}$, where Bethe vectors in the Gaudin model were
produced. Solutions of the Bethe ansatz equations associated with
the general orthosymplectic Lie superalgebras were investigated
in \cite{lm:ba} via a reproduction procedure.

Note that the evaluation homomorphism \eqref{evalrgm}
can be modified by using the additional parameter $\mu\in\g^*$
such that $\mu$ vanishes on the odd elements of $\g$. Equivalently,
$\mu$ can be regarded as an element of $\g_{\bar 0}^*$ with $\g_{\bar 0}=\oa_M\oplus\spa_{2n}$.
The modified homomorphism takes the form
\beql{evalrgmod}
\U\big(t^{-1}\g[t^{-1}]\big)\to \U(\g),
\qquad F_{ij}[-r]\mapsto F_{ij}\tss z^{-r}+\de_{r1}\ts\mu(F_{ij}),\quad r>0.
\eeq
This leads to a more general theory of Gaudin algebras as originally considered in \cite{r:si}
and further developed in \cite{fft:gm}. It would be interesting to extend the work
\cite{mm:eb} to the orthosymplectic case by calculating the eigenvalues of Bethe vectors;
see also \cite[Ch.~14]{m:so}.

\subsection{Quantum Mishchenko--Fomenko subalgebras}
\label{subsec:qmf}

The classical construction of Poisson commutative subalgebras of the symmetric
algebra $\Sr(\g)$ of a Lie algebra $\g$, originated in \cite{mf:ee}, can be naturally
extended to Lie superalgebras. Regarding the elements $F_{ij}$ of the
orthosymplectic Lie superalgebra $\g=\osp_{M|2n}$ as generators
of the symmetric superalgebra $\Sr(\g)$, suppose that $P\in\Sr(\g)^{\g}$
is a $\g$-invariant of $\Sr(\g)$ under the adjoint action.
Take any element $\mu\in\g^*$ which
vanishes on the odd elements of $\g$. Regarding $P$
as a polynomial in the $F_{ij}$, use a `shift of argument'
to replace the variables by
$
F_{ij}\mapsto F_{ij}+t\tss \mu(F_{ij}),
$
where $t$ is a variable. After this replacement, the new polynomial expands
as a polynomial in $t$,
\beql{polte}
P^{}_{(0)}+P^{}_{(1)}\tss t+\dots+P^{}_{(k)}\tss t^k,
\eeq
thus defining
elements $P^{}_{(i)}\in \Sr(\g)$ associated with $P$ and $\mu$.
The {\em (classical) Mishchenko--Fomenko subalgebra}
$\overline\Ac_{\mu}$ of $\Sr(\g)$
is generated by all elements $P^{}_{(i)}$
associated with all $\g$-invariants $P\in \Sr(\g)^{\g}$.
The key property of the subalgebra $\overline\Ac_{\mu}$, which is verified in the same way
as in the Lie algebra case, is that $\overline\Ac_{\mu}$ is Poisson commutative
with respect to the Lie--Poisson super-bracket on $\Sr(\g)$; see e.g. \cite{r:si}
and \cite[Sec.~9.1]{m:so}. One can therefore extend {\em Vinberg's quantization problem}
\cite{v:sc}
by asking whether it is possible to
construct a commutative
subalgebra $\Ac_{\mu}$ of $\U(\g)$ that `quantizes' $\overline\Ac_{\mu}$.
Here we equip $\U(\g)$ with the canonical filtration and
regard $\U(\g)$ as a `quantization'
of the Poisson superalgebra $\Sr(\g)$,
in the sense that the graded superalgebra $\gr\U(\g)$
is isomorphic to $\Sr(\g)$. The required
quantization property for the Mishchenko--Fomenko subalgebra
then reads $\gr\Ac_{\mu}=\overline\Ac_{\mu}$.

A universal solution of Vinberg's quantization problem
based on the use of the Feigin--Frenkel centre $\z(\wh\g)$ was
proposed in \cite{r:si} in the Lie algebras case; see also \cite{fft:gm} and
\cite[Ch.~9]{m:so} for further results and more references.
Namely, the subalgebra $\Ac_{\mu}$ of $\U(\g)$ is defined as a homomorphic image
of $\z(\wh\g)$, and conjecturally, the property $\gr\Ac_{\mu}=\overline\Ac_{\mu}$
holds; see \cite[Conjecture~1]{fft:gm}. The conjecture was proved
for regular $\mu$ in \cite{fft:gm}, while a proof for type $C$ and
a new proof in type $A$
for all $\mu$ was given in \cite{my:qn}.

We use the same approach for the orthosymplectic Lie superalgebra $\g$ and
define the subalgebra $\Ac_{\mu}\subset\U(\g)$ as the image
of $\z(\wh\g)$ with respect to the
homomorphism
\eqref{evalrgmod}. The image is easily seen not to depend on $z$.
Explicit formulas for elements of $\Ac_{\mu}$ then
follow from the results of \cite{mn:ss}; see the comments
following \eqref{poltaexpa}.
Introduce the matrix $\mu=[\mu(\hF_{ij})]$.

\bpr\label{prop:qmf}
The coefficients of all polynomials in $z^{-1}$ given by
\ben
\ga_k(M-2n)\ts\str^{}_{1,\dots,k}\tss S^{(k)} \big({-}\di_z+\mu_1+\hF_{1}z^{-1}\big)\dots
\big({-}\di_z+\mu_k+\hF_{k}z^{-1}\big)\tss 1
\een
with positive integer values of $k$
belong to the commutative superalgebra $\Ac_{\mu}$.
\qed
\epr

Note that it is implied by \cite[Conjecture~2.3]{mn:ss} that
the elements defined in Proposition~\ref{prop:qmf} generate the superalgebra
$\Ac_{\mu}$ in the case of odd $M$. If $M$ is even, then additional elements arising
from the conjectural super-Pfaffian Segal--Sugawara vector $\Phi'$ should be necessary
to generate $\Ac_{\mu}$, as stated in Conjecture~\ref{conj:even}.
In all cases, we expect the following to hold.

\bcj\label{conj:mf}
The superalgebras $\Ac_{\mu}$ solve Vinberg's quantization problem:
$\gr\Ac_{\mu}=\overline\Ac_{\mu}$ for all $\mu$.
\ecj

\section{Proof of Theorem~\ref{thm:hch}}
\label{sec:pt}

Our proof will rely on the connection between the Yangian characters and
Harish-Chandra images of the Segal--Sugawara vectors that
was already used in \cite{mm:yc} to prove the counterparts
of Theorem~\ref{thm:hch} for the orthogonal and symplectic Lie algebras; see also
\cite[Ch.~11 and 13]{m:so}. The {\em completed dual Yangian} $\wh\Y^+(\osp_{M|2n})$
possesses a filtration such that the associated graded algebra $\gr\wh\Y^+(\osp_{M|2n})$
is isomorphic to the universal enveloping algebra $\U(t^{-1}\osp_{M|2n}[t^{-1}])$.
The isomorphism turns out to respect the Harish-Chandra homomorphisms
for the Yangian and the universal enveloping algebra. The Yangian characters
are the Harish-Chandra images of certain formal series with coefficients in the Yangian, and they
can be evaluated with the use of $R$-matrix calculations by employing special bases
of the space of orthosymplectic harmonic polynomials. We
find suitable linear combinations of
the Yangian characters and calculate their top degree terms as elements of the associated
graded algebra, thus producing the Harish-Chandra images of
the Segal--Sugawara vectors $\Phi_k$.
We start by calculating certain Yangian characters for the
orthosymplectic Yangian to follow by their dual counterparts.

\subsection{Yangian characters}
\label{subsec:yc}

The {\em orthosymplectic Yangian} $\Y(\g)$ for $\g=\osp_{M|2n}$
was originally introduced in \cite{aacfr:rp};
we will recall the definitions following the settings of \cite{mr:rb}.
The $R$-{\em matrix} associated with $\g$ is the
rational function in $u$ given by
\beql{zamolr}
R(u)=1-\frac{P}{u}+\frac{Q}{u-\ka},\qquad \ka=\frac{M}{2}-n-1,
\eeq
where $P$ and $Q$ are defined in Sec.~\ref{subsec:def}.
The {\em extended Yangian}
$\X(\g)$ is
the $\ZZ_2$-graded algebra with generators
$t_{ij}^{(r)}$ of parity $\bi+\bj\mod 2$, where $1\leqslant i,j\leqslant 1'$ and $r=1,2,\dots$,
satisfying the following defining relations.
Introduce the formal series
\beql{tiju}
t_{ij}(u)=\de_{ij}+\sum_{r=1}^{\infty}t_{ij}^{(r)}\ts u^{-r}
\in\X(\g)[[u^{-1}]]
\eeq
and combine them into the matrix $T(u)=[t_{ij}(u)]$.
Consider the elements of the tensor product superalgebra
$\End\CC^{M|2n}\ot\End\CC^{M|2n}\ot \X(\g)[[u^{-1}]]$ given by
\beql{tonettwo}
T_1(u)=\sum_{i,j=1}^{1'} e_{ij}\ot 1\ot t_{ij}(u)(-1)^{\bi\tss\bj+\bj}\fand
T_2(u)=\sum_{i,j=1}^{1'} 1\ot e_{ij}\ot t_{ij}(u)(-1)^{\bi\tss\bj+\bj}.
\eeq
The defining relations for $\X(\g)$ take
the form of the $RTT$-{\em relation}
\beql{RTT}
R(u-v)\ts T_1(u)\ts T_2(v)=T_2(v)\ts T_1(u)\ts R(u-v).
\eeq
The {\em Yangian} $\Y(\g)$
is defined as the quotient
of $\X(\g)$
by the relation $T(u+\ka)^{\tss t}\ts T(u)=1$.

We will identify the universal enveloping algebra $\U(\g)$
with a subalgebra of $\X(\g)$ via the embedding
\beql{embla}
F_{ij}\mapsto \frac12\big(t_{ij}^{(1)}-t_{j'i'}^{(1)}(-1)^{\bi\bj+\bj}\ts\ta_i\ta_j\big)(-1)^{\bi}.
\eeq
The image of $F_{ij}$ can also be written as $t_{ij}^{(1)}(-1)^{\bi}+\de_{ij}\tss \ze$
for a central element $\zeta$. Note that $\zeta=0$ in $\Y(\g)$.
Then $\X(\g)$ can be regarded as a $\g$-module
with the adjoint action.
Recall the triangular decomposition \eqref{tridec} and
denote by $\X(\g)^{\h}$ the
subalgebra of $\h$-invariants under this action,
\ben
\X(\g)^{\h}=\{y\in \X(\g)\ |\ [F_{ii},y]=0\quad\text{for}\quad i=1,\dots,n+m\}.
\een
Consider the left ideal $\Ir$ of the algebra $\X(\g)$ generated by
all elements $t^{(r)}_{ij}$ with the conditions
$1\leqslant i<j\leqslant 1'$ and $r\geqslant 1$.
The intersection
$\X(\g)^{\h}\cap \Ir$ is a two-sided ideal of $\X(\g)^{\h}$. Moreover,
the quotient of $\X(\g)^{\h}$ by this ideal is isomorphic
to a commutative algebra generated by the images
of the elements $t_{ii}^{(r)}$ with $i=1,\dots,1'$ and $r\geqslant 1$
in the quotient. We will use the notation $\la^{(r)}_i$
for the image of $t_{ii}^{(r)}$. Thus, we get the Harish-Chandra
homomorphism
\beql{yhchbcd}
\X(\g)^{\h}\to\text{\quad\rm polynomials in\quad}\la^{(r)}_i,\quad i=1,\dots,1',\ r\geqslant 1.
\eeq
We combine the elements $\la^{(r)}_i$ into the formal series
\beql{laiugn}
\la_i(u)=1+\sum_{r=1}^{\infty} \la^{(r)}_i\ts u^{-r},\qquad i=1,\dots,1',
\eeq
which can be understood as the images of the series
$t_{ii}(u)$ under the homomorphism \eqref{yhchbcd}.

Recall the action
of the Brauer algebra $\Bc_k(\om)$ with $\om=M-2n$
on the superspace $(\CC^{M|2n})^{\ot k}$
which was already used in \cite{mn:ss}. The algebra $\Bc_k(\om)$ has a basis of diagrams, where each diagram consists of two horizontal rows of $k$ nodes, and $k$ strings connecting the nodes pairwise. The product $xy$ of two diagrams $x$ and $y$ is computed by concatenation; we draw $y$ directly above $x$, connect the strings at the nodes in the middle, remove the middle nodes, 
and replace each loop formed by a factor of $\om$. It is well-known that
the algebra $\Bc_k(\om)$ is generated by the set of diagrams
\begin{align}
    s_{ab} &= \, \begin{tikzpicture}[baseline={([yshift=4.5mm]current bounding box.south)},scale=0.5]
{
\foreach \x in {0,...,7}{
\filldraw (\x,0) circle (0.09);
\filldraw (\x,2) circle (0.09);
}
\draw (0,0)--(0,2);
\draw (1,0)--(1,2);
\draw (3,0)--(3,2);
\draw (4,0)--(4,2);
\draw (6,0)--(6,2);
\draw (7,0)--(7,2);
\draw (2,0) to[out=90,in=270](5,2);
\draw (5,0) to[out=90,in=270](2,2);
\foreach \y in {0.2,0,-0.2}{
\node at (0.5+\y,1) {.};
}
\foreach \y in {0.2,0,-0.2}{
\node at (3.5+\y,0.5) {.};
}
\foreach \y in {0.2,0,-0.2}{
\node at (6.5+\y,1) {.};
}
\draw (0,2) node[font=\scriptsize,anchor=south]{$1$};
\draw (2,2) node[font=\scriptsize,anchor=south]{$a$};
\draw (5,2) node[font=\scriptsize,anchor=south]{$b$};
\draw (7,2) node[font=\scriptsize,anchor=south]{$k$};
}
\end{tikzpicture}\; ,
    &\epsilon_{ab} &= \, \begin{tikzpicture}[baseline={([yshift=4.5mm]current bounding box.south)},scale=0.5]
{
\foreach \x in {0,...,7}{
\filldraw (\x,0) circle (0.09);
\filldraw (\x,2) circle (0.09);
}
\draw (0,0)--(0,2);
\draw (1,0)--(1,2);
\draw (3,0)--(3,2);
\draw (4,0)--(4,2);
\draw (6,0)--(6,2);
\draw (7,0)--(7,2);
\draw (2,0) ..controls (2,0.9) and (5,0.9).. (5,0);
\draw (2,2) ..controls (2,1.1) and (5,1.1).. (5,2);
\foreach \y in {0.2,0,-0.2}{
\node at (0.5+\y,1) {.};
}
\foreach \y in {0.2,0,-0.2}{
\node at (3.5+\y,1) {.};
}
\foreach \y in {0.2,0,-0.2}{
\node at (6.5+\y,1) {.};
}
\draw (0,2) node[font=\scriptsize,anchor=south]{$1$};
\draw (2,2) node[font=\scriptsize,anchor=south]{$a$};
\draw (5,2) node[font=\scriptsize,anchor=south]{$b$};
\draw (7,2) node[font=\scriptsize,anchor=south]{$k$};
}
\end{tikzpicture}\;
\non
\end{align}
with $1\leqslant a<b\leqslant k$.
For each $1 \leqslant l \leqslant k$, we define the
{\em partial transposition} $t_l:\Bc_k(\om) \to \Bc_k(\om)$
as the linear map taking each basis diagram $d$ to the diagram obtained
by swapping the nodes numbered $l$ on the top and on the bottom lines of $d$
while keeping the strings attached to them. For example, in $\Bc_6(\om)$, we have
\begin{align}
    d &= \begin{tikzpicture}[baseline={([yshift=4.5mm]current bounding box.south)},scale=0.5]
{
\foreach \x in {0,...,5}{
\filldraw (\x,0) circle (0.09);
\filldraw (\x,2) circle (0.09);
}
\draw (0,0) to[out=90,in=270] (2,2);
\draw (1,0) to[out=90,in=270] (0,2);
\draw (2,0) ..controls (2,0.7) and (5,0.7).. (5,0);
\draw (3,0) to[out=90,in=270] (5,2);
\draw (4,0) to[out=90,in=270] (3,2);
\draw (1,2) to[out=270,in=270] (4,2);
}
\end{tikzpicture}\; ,
    &d^{t_5} &= \begin{tikzpicture}[baseline={([yshift=4.5mm]current bounding box.south)},scale=0.5]
{
\foreach \x in {0,...,5}{
\filldraw (\x,0) circle (0.09);
\filldraw (\x,2) circle (0.09);
}
\draw (0,0) to[out=90,in=270] (2,2);
\draw (1,0) to[out=90,in=270] (0,2);
\draw (2,0) ..controls (2,0.7) and (5,0.7).. (5,0);
\draw (3,0) to[out=90,in=270] (5,2);
\draw (4,0) to[out=90,in=270] (1,2);
\draw (3,2) arc(180:360:0.5);
}
\end{tikzpicture}\; .
\non
\end{align}
There is a homomorphism
\beql{brahom}
\Bc_{k}(\om) \to \End(\CC^{M|2n})^{\otimes k}
\eeq
which is defined by $s_{ab} \mapsto P_{ab}$ and $\epsilon_{ab} \mapsto Q_{ab}$,
where the operators $P_{ab}$ and $Q_{ab}$ are given by
\begin{align}
P_{ab} &= \sum_{i,j=1}^{M+2n} 1^{\otimes (a-1)} \otimes e_{ij} \otimes 1^{\otimes (b-a-1)} \otimes e_{ji} \otimes 1^{\otimes (k-b)} (-1)^{\bj},
\label{pab}\\
Q_{ab} &= \sum_{i,j=1}^{M+2n} 1^{\otimes (a-1)} \otimes e_{ij} \otimes 1^{\otimes (b-a-1)} \otimes e_{i'j'} \otimes 1^{\otimes (k-b)} (-1)^{\bi\bj} \ta_i\ta_j.
\label{qab}
\end{align}

Following \cite{mn:ss} and \cite[Sec.~1.2]{m:so}, consider the symmetrizer $s^{(k)}$ in the Brauer algebra
which is an idempotent associated with the trivial representation of $\Bc_k(\om)$.
Its image $S^{(k)}$ under the action of $\Bc_k(M-2n)$ on the tensor product
superspace $(\CC^{M|2n})^{\ot k}$ is well-defined for those values
of the parameters where the denominators occurring in the formula for $s^{(k)}$
do not vanish.
We will identify $S^{(k)}$ with
the element
$S^{(k)}\ot 1$ of the algebra
\beql{tenprky}
\underbrace{\End\CC^{M|2n}\ot\dots\ot\End\CC^{M|2n}}_k{}\ot\X(\g).
\eeq
Explicitly, it can be given by the multiplicative
formula
\beql{symo}
S^{(k)}=\frac{1}{k!}
\prod_{1\leqslant a<b\leqslant k}
\Big(1+\frac{P_{a\tss b}}{b-a}-\frac{Q_{a\tss b}}
{\ka+b-a}\Big),
\eeq
where
the products are taken in the lexicographic order
on the pairs $(a,b)$.
Note that
the products involve evaluated $R$-matrices \eqref{zamolr},
\beql{rmapr}
S^{(k)}=\frac{1}{k!}
\prod_{1\leqslant a<b\leqslant k} R_{a\tss b}(u_a-u_b),
\eeq
where $u_a=u+a-1$ for $a=1,\dots,k$.

If $M=2m+1$, then the expression in \eqref{symo} is defined for all values
of the parameters. If $M=2m$, then we will assume that $m>n$ for the expression to make sense.

Introduce the formal series $\TT^{(k)}(u)$ with coefficients in
the extended Yangian $\X(\g)$ by the formula
\beql{betheb}
\TT^{(k)}(u)=\str_{1,\dots,k}\ts S^{(k)}T_1(u)\ts T_2(u+1)\dots T_k(u+k-1)
\eeq
with the supertrace taken over all $k$ copies of $\End\CC^{M|2n}$ in \eqref{tenprky},
where we extend notation \eqref{tonettwo} to copies of the matrix $T(u)$
in \eqref{tenprky}.
All coefficients of this series
belong to $\X(\g)^{\h}$.

The following theorem provides the formulas for the required Yangian characters
in the case of odd $M$; cf. \cite[Sec.~11.1]{m:so}.

\bth\label{thm:hchb}
If $M=2m+1$ then
the image of the series $\TT^{(k)}(u)$ under the homomorphism \eqref{yhchbcd}
is given by
\ben
\TT^{(k)}(u)\mapsto\sum_{i_1<\dots<i_s\leqslant n<i_{s+1}\leqslant\dots\leqslant i_p<n'\leqslant i_{p+1}<\dots<i_k}
\la_{i_1}(u)\ts\la_{i_2}(u+1)\dots \la_{i_k}(u+k-1)(-1)^{\bi_1+\dots+\bi_k},
\een
assuming varying values of the parameters
$s$ and $p$ in the summation,
with the condition that $n+m+1$ occurs among the summation indices $i_1,\dots,i_k$
at most once.
\eth

\bpf
The subspace of {\em harmonic tensors}
in $(\CC^{M|2n})^{\ot k}$ is spanned
by the tensors $v$ with the property $Q_{a\tss b}\ts v=0$ for all $1\leqslant a<b\leqslant k$.
By the properties of the symmetrizer $s^{(k)}$ described in \cite[Sec.~1.2]{m:so},
the operator $S^{(k)}$ projects $(\CC^{M|2n})^{\ot k}$ to a subspace of {\em supersymmetric}
harmonic tensors $\Hc_k=S^{(k)}(\CC^{M|2n})^{\ot k}$. This subspace was considered in
\cite[Cor.~1]{ds:sh}, where its dimension was calculated, and in \cite{c:os},
where it was described as an $\osp_{M|2n}$-module. For $M\geqslant 1$ we have
\beql{dimhk}
\dim \Hc_k=\sum_{i=0}^{\min(k,2n)}\binom{2n}{i}\binom{M+k-i-1}{M-1}
-\sum_{i=0}^{\min(k-2,2n)}\binom{2n}{i}\binom{M+k-i-3}{M-1}.
\eeq

We will keep notation
$H^{(k)}$ for the image of the element \eqref{ha} under the action of the Brauer algebra
on the tensor product superspace.
Denote by $\Pc_{M|2n}$
the supercommutative superalgebra with generators $z_1,z_2,\dots,z_{1'}$,
where the parity of $z_i$ is defined as $\bi\mod 2$.
Identify the image of the tensor product superspace under the operator
$H^{(k)}$ with the subspace $\Pc_{M|2n}^{\tss k}$
of homogeneous elements of degree $k$,
\beql{isomhmcm}
H^{(k)}(\CC^{M|2n})^{\ot k}\cong \Pc_{M|2n}^{\tss k}
\eeq
by setting
\beql{symsubo}
H^{(k)}(e_{i_1}\ot\dots\ot e_{i_k})=z_{i_1}\dots z_{i_k}.
\eeq
As with the non-super case considered in \cite[Ch.~2]{m:so}, we may regard
$S^{(k)}$ as an operator on $\Pc_{M|2n}^{\tss k}$. Moreover, its action
coincides with that of the {\em extremal projector} $p$ associated with the Lie algebra $\sll_2$.
According to a particular case of the {\em Howe duality},
the standard basis elements $e,f,h$ of $\sll_2$ act on $\Pc_{M|2n}$ by the rules
\beql{sltwoa}
\bal
e&= \sum_{i=1}^n \di_{i}\tss
\di_{\tss i^{\tss\prime}}-\frac12\ts\sum_{i=n+1}^{(n+1)'} \di_{i}\tss
\di_{\tss i^{\tss\prime}},\\
f&= -\sum_{i=1}^n z_{i}\tss
z_{\tss i^{\tss\prime}}+\frac12\ts\sum_{i=n+1}^{(n+1)'} z_i\tss z_{\tss {i\tss}'},\\
h&= -\frac{M}{2}+n-\sum_{i=1}^{1'}
z_i\ts\di_i,
\eal
\eeq
where $\di_i$ denotes the partial derivative over $z_i$. Hence, using
the correspondence \eqref{symsubo} we can identify
the space $\Hc_k$ of supersymmetric harmonic tensors with
the subspace of {\em harmonic polynomials} of degree $k$ in $z_1,\dots,z_{1'}$.
This subspace is defined as the kernel
of the {\em Laplace operator} $e$ so that
\beql{harmo}
\Hc_k\cong \{P\in \Pc_{M|2n}^{\tss k}\ts |\ts e\tss P=0\}.
\eeq

As a next step, we produce a basis $P_I$ of $\Hc_k$ which will be parameterized
by the multisets $I=\{i_1,\dots,i_k\}$ such that
\beql{msi}
1\leqslant i_1<\dots<i_s\leqslant n<i_{s+1}\leqslant\dots\leqslant
i_p<n'\leqslant i_{p+1}<\dots<i_k\leqslant 1'
\eeq
for all possible $s$ and $p$, with the condition that $n+m+1$ occurs
in the multiset at most once.
Given such a multiset $I$, consider the corresponding monomial
$z_I=z_{i_1}\dots z_{i_k}$ and write it as $z_I=z_{n+m+1}^{\de}z_{\bar I}$,
where $\de=1$ if $n+m+1$ occurs in the multiset and $\de=0$ otherwise,
with $\bar I=I\setminus \{n+m+1\}$ and $\bar I=I$, respectively.

Now split the Laplace operator by setting
\ben
e=\bar{e}+\De,\qquad \bar{e}=\sum_{i=1}^n \di_{i}\tss
\di_{\tss i^{\tss\prime}}-\sum_{i=n+1}^{n+m} \di_{i}\di_{\tss i^{\tss\prime}},\qquad
\De=-\frac12\ts\di_{n+m+1}^2.
\een
Define elements $P_I\in\Pc_{M|2n}^{\tss k}$ by
\beql{basb}
P_I=\sum_{s\geqslant 0}\frac{2^s}{(\de+2s)!}\ts z_{n+m+1}^{\de+2s}{\bar{e}}^{\ts s}\ts z_{\bar I}.
\eeq
Let us verify that the polynomials $P_I$ are all harmonic and form a basis of $\Hc_k$.
Since $\De z_{I}=0$,
by applying the Laplace operator, we get
\ben
e\ts P_I=(\bar{e}+\De)\ts P_I
=\sum_{s\geqslant 0}\frac{2^s}{(\de+2s)!}\ts z_{n+m+1}^{\de+2s}\ts{\bar{e}}^{\ts s+1}\ts z_{\bar I}
-\frac12\ts \sum_{s\geqslant 1}\frac{2^s}{(\de+2s-2)!}\ts z_{n+m+1}^{\de+2s-2}
\ts{\bar{e}}^{\ts s}\ts z_{\bar I}=0.
\een
The polynomials $P_I$ are linearly independent because each $P_I$ contains a unique monomial
(the {\em leading
monomial}) where the variable $z_{n+m+1}$ occurs with
the power not exceeding $1$. Finally, the number of possible multisets $I$ equals
\ben
\sum_{i=0}^{\min(k,2n)}\binom{2n}{i}\Bigg(\binom{M+k-i-2}{k-i}+\binom{M+k-i-3}{k-i-1}\Bigg)
\een
which coincides with $\dim \Hc_k$ in \eqref{dimhk}.

To complete the proof of the theorem, note that
by \eqref{RTT} and \eqref{rmapr},
the product
occurring in \eqref{betheb} can be written as
\beql{stt}
S^{(k)}T_1(u)\dots T_k(u+k-1)=
T_k(u+k-1)\dots T_1(u)\ts S^{(k)}.
\eeq
Hence the product on
each side can be regarded as
an operator on $(\CC^{M|2n})^{\ot k}$
such that the subspace $\Hc_k$ is invariant under this operator.
Therefore, the supertrace in \eqref{betheb} can be
calculated over the subspace $\Hc_k$. To do this,
fix a basis vector $P_I\in \Hc_k$ and use the correspondence \eqref{symsubo}
to regard $P_I$ as an element of tensor product space $(\CC^{M|2n})^{\ot k}$,
\beql{expapi}
P_I=\sum_J c_{J}\ts e_{j_1}\ot\dots\ot e_{j_k}, \qquad J=(j_1,\dots,j_k).
\eeq
Since $S^{(k)}P_I=P_I$, the application of the operator appearing on the right hand side
of \eqref{stt} to $P_I$ yields
\begin{multline}
\sum_J c_{J}\ts T_k(u+k-1)\dots T_1(u)(e_{j_1}\ot\dots\ot e_{j_k})\\[0.3em]
{}=\sum_J c_{J}\ts\sum_{h_1,\dots,h_k}(e_{h_1}\ot\dots\ot e_{h_k})\ts
t_{h_kj_k}(u+k-1)\dots t_{h_1j_1}(u)(-1)^{\sum_{a<b}(\bh_a+\bj_a)\bj_b}.
\non
\end{multline}
This equals a linear combination of the basis vectors of $\Hc_k$, and we want to
evaluate the coefficient of the basis vector $P_I$.
The coefficient of $P_I$ in the linear combination is uniquely determined
by the coefficient of the tensor $e_I=e_{i_1}\ot\dots\ot e_{i_k}$.
It follows from
\eqref{basb} that if a tensor of the form
$e_{j_1}\ot\dots\ot e_{j_k}$ corresponds to
a non-leading monomial occurring in the expansion \eqref{expapi},
then the matrix element
\ben
t_{i_kj_k}(u+k-1)\dots t_{i_1j_1}(u)
\een
vanishes under the homomorphism \eqref{yhchbcd}.
Indeed, if $i_1\leqslant n$, and the coefficient $c_J$ in \eqref{expapi}
is nonzero, then $j_1\geqslant i_1$. Hence,
if we suppose that the matrix element does not vanish
under the homomorphism \eqref{yhchbcd}, then we must have $j_1=i_1$.
Similarly, if $i_1<i_2\leqslant n$, then $j_2=i_2$, etc., concluding that
$j_a=i_a$ for all $a=1,\dots,s$ in \eqref{msi}.
Furthermore,
the multiplicity of the index $n+1$ in the multiset $\{j_1,\dots,j_k\}$
does not exceed its multiplicity in the multiset $I$, whereas
the coefficients of the series $t_{n+1, j_c}(u+c-1)$
belong to the left ideal $\Ir$ of $\X(\g)$if $n+1< j_c$. Hence,
if we suppose that the matrix element does not vanish
under the homomorphism \eqref{yhchbcd}, then
the multiplicities of $n+1$ in the multisets $\{j_1,\dots,j_k\}$
and $I$ must coincide. By extending this argument to $n+2,\dots, n+m$, we can conclude
that the multiplicities of each of these indices in the multisets $\{j_1,\dots,j_k\}$
and $I$ must coincide.

By the formula for the polynomial $P_I$ as defined in \eqref{basb},
the same conclusion can then be made regarding the multiplicities of each of the
indices $(n+m)',\dots,1'$. This is clear from the form of the operator $\bar{e}$,
as the products $\di_{i}\di_{\tss i^{\tss\prime}}$ act on the polynomials by
simultaneously decreasing the degrees with respect to $z_i$ and $z_{i'}$
for $i=1,\dots,n+m$.

Thus,
a nonzero contribution to the
image of the diagonal matrix element of the operator corresponding
to $P_I$ comes only from
the term  $t_{i_ki_k}(u+k-1)\dots t_{i_1i_1}(u)$.
Calculating the supertrace,
sum over all basis vectors \eqref{basb} taking the parity into account
to get the desired
formula for the image of the element \eqref{betheb}.
\epf

Now we prove a counterpart of Theorem~\ref{thm:hchb} for even
values of $M$.

\bth\label{thm:hchd}
If $M=2m$ and $m>n$, then
the image of the series $\TT^{(k)}(u)$ under the homomorphism \eqref{yhchbcd}
is given by
\ben
\TT^{(k)}(u)\mapsto\sum_{i_1<\dots<i_s\leqslant n<i_{s+1}\leqslant\dots\leqslant i_p<n'\leqslant i_{p+1}<\dots<i_k}
\la_{i_1}(u)\ts\la_{i_2}(u+1)\dots \la_{i_k}(u+k-1)(-1)^{\bi_1+\dots+\bi_k},
\een
assuming varying values of the parameters
$s$ and $p$ in the summation,
with the condition that $n+m$ and $(n+m)'$ do not occur
simultaneously
among the summation indices $i_1,\dots,i_k$.
\eth

\bpf
Since the dimension formula \eqref{dimhk} for the harmonic polynomials
is valid for even values $M\geqslant 2$, and the Lie algebra
$\sll_2$ acts by the same formulas \eqref{sltwoa},
the starting arguments in the proof of Theorem~\ref{thm:hchb} apply
in the case of even $M$ as well. However, we will need to
modify the construction of the basis of the space of harmonic polynomials
occurring in \eqref{harmo}.
Our basis $P_I$ of $\Hc_k$ will now be parameterized
by the multisets $I=\{i_1,\dots,i_k\}$ such that
\beql{msid}
1\leqslant i_1<\dots<i_s\leqslant n<i_{s+1}\leqslant\dots\leqslant
i_p<n'\leqslant i_{p+1}<\dots<i_k\leqslant 1'
\eeq
for all possible $s$ and $p$, with the condition that $n+m$ and $(n+m)'$ do not occur
simultaneously
in the multiset. Given such a multiset, introduce the corresponding monomial
$z_I=z_{i_1}\dots z_{i_k}$.
This time we split the Laplace operator by
\ben
e=\bar{e}+\De,\qquad \bar{e}=\sum_{i=1}^n \di_{i}\tss
\di_{\tss i^{\tss\prime}}-\sum_{i=n+1}^{n+m-1} \di_{i}\di_{\tss i^{\tss\prime}},\qquad
\De=-\di_{n+m}\di_{\tss (n+m)^{\tss\prime}}.
\een
Define elements $P_I\in\Pc_{M|2n}^{\tss k}$ by
\beql{basd}
P_I=\sum_{s\geqslant 0}\frac{z_{n+m}^s}{(s+l_{n+m})!}
\frac{z_{(n+m)'}^s}{(s+l_{(n+m)'})!}{\bar{e}}^{\ts s}\ts z_{I},
\eeq
where $l_{n+m}$ and $l_{(n+m)'}$ are the respective multiplicities of $n+m$ and $(n+m)'$
in the multiset $I$; at least one of them is zero.

As with the polynomials \eqref{basb}, we use the property $\De z_{I}=0$ to verify
that the polynomials $P_I$ are all harmonic.
They are linearly independent because each $P_I$ contains a unique monomial
(the {\em leading
monomial}) which does not contain at least one of $z_{n+m}$ or $z_{(n+m)'}$.
The number of multisets \eqref{msid} equals
\ben
\sum_{i=0}^{\min(k,2n)}\binom{2n}{i}\Bigg(2\binom{M+k-i-2}{k-i}-\binom{M+k-i-3}{k-i}\Bigg)
\een
which coincides with $\dim \Hc_k$ in \eqref{dimhk}, thus proving that the $P_I$
form a basis of $\Hc_k$.

The final argument is the same as for Theorem~\ref{thm:hchb}; we show that
a nonzero contribution to the Harish-Chandra
image of the diagonal matrix element of the operator on the basis vector
$P_I$ comes only from
the term  $t_{i_ki_k}(u+k-1)\dots t_{i_1i_1}(u)$.
\epf

We will also need dual versions of
Theorems~\ref{thm:hchb} and \ref{thm:hchd},
where the Harish-Chandra homomorphism is defined with respect to
the left ideal $\Ir^{\op}$ instead of $\Ir$. Namely,
the left ideal $\Ir^{\op}$ of the algebra $\X(\g)$ is generated by
all elements $t^{(r)}_{ij}$ with the conditions
$1'\geqslant i>j\geqslant 1$ and $r\geqslant 1$.
Similar to \eqref{yhchbcd}, we get the Harish-Chandra
homomorphism
\beql{yhchbcdop}
\X(\g)^{\h}\to\text{\quad\rm polynomials in\quad}\la^{(r)}_i,\quad i=1,\dots,1',\ r\geqslant 1,
\eeq
defined as the projection modulo the two-sided ideal $\X(\g)^{\h}\cap \Ir^{\op}$,
where we use the same notation for
the formal series \eqref{laiugn},
understood as the images of the respective series
$t_{ii}(u)$ under the homomorphism \eqref{yhchbcdop}.

\bco\label{cor:hchbdy}
If $M=2m+1$ then
the image of the series $\TT^{(k)}(u)$ under the homomorphism \eqref{yhchbcdop}
is given by
\ben
\TT^{(k)}(u)\mapsto\sum_{i_1>\dots>i_s\geqslant n'>i_{s+1}\geqslant\dots\geqslant i_p>n\geqslant i_{p+1}>\dots>i_k}
\la_{i_1}(u)\ts\la_{i_2}(u+1)\dots \la_{i_k}(u+k-1)(-1)^{\bi_1+\dots+\bi_k},
\een
assuming varying values of the parameters $s$ and $p$ in the summation,
with the condition that $n+m+1$ occurs among the summation indices $i_1,\dots,i_k$
at most once.
\eco

\bco\label{cor:hchdy}
If $M=2m$ and $m>n$ then
the image of the series $\TT^{(k)}(u)$ under the homomorphism \eqref{yhchbcdop}
is given by
\ben
\TT^{(k)}(u)\mapsto\sum_{i_1>\dots>i_s\geqslant n'>i_{s+1}\geqslant\dots\geqslant i_p>n\geqslant i_{p+1}>\dots>i_k}
\la_{i_1}(u)\ts\la_{i_2}(u+1)\dots \la_{i_k}(u+k-1)(-1)^{\bi_1+\dots+\bi_k},
\een
assuming varying values of the parameters $s$ and $p$ in the summation,
with the condition that $n+m$ and $(n+m)'$ do not occur
simultaneously
among the summation indices $i_1,\dots,i_k$.
\eco

\bpf[Proof of Corollaries~\ref{cor:hchbdy} and \ref{cor:hchdy}]
The corollaries
follow from Theorems~\ref{thm:hchb} and \ref{thm:hchd}, respectively,
with the use of the automorphism of the extended Yangian $\X(\g)$
defined by
\beql{autoop}
t_{ij}(u)\mapsto t_{i'j'}(u)\ts\ta_i\ta_j.
\eeq
Obviously, the automorphism takes the left ideal $\Ir$ to $\Ir^{\op}$
and $t_{ii}(u)$ to $t_{i'i'}(u)$, so that the desired
formulas for the Harish-Chandra images follows by taking
the composition of the homomorphism \eqref{yhchbcd}
with the automorphism \eqref{autoop}.
\epf

\subsection{From dual Yangians to vacuum modules}
\label{subsec:dy}

Define the
{\em extended dual Yangian}
$\X^+(\g)$ for $\g=\osp_{M|2n}$
as the $\ZZ_2$-graded algebra with generators
$t_{ij}^{(-r)}$ of parity $\bi+\bj\mod 2$,
where $1\leqslant i,j\leqslant 1'$ and $r=1,2,\dots$,
subject to the defining relations written in a matrix form as follows.
Combine the generators into the formal power series
\ben
t^+_{ij}(u)=\de_{ij}-\sum_{r=1}^{\infty}t_{ij}^{(-r)}\ts u^{r-1}
\in\X^+(\g)[[u]]
\een
and introduce the matrix
\ben
T^+(u)=\sum_{i,j=1}^{1'} e_{ij}\ot t^+_{ij}(u)(-1)^{\bi\tss\bj+\bj}
\in \End\CC^{M|2n}\ot \X^+(\g)[[u]].
\een
The defining relations are
\beql{RTTdualbcd}
R_{12}(u-v)\ts T^+_1(u)\ts T^+_2(v)=T^+_2(v)\ts T^+_1(u)\ts R_{12}(u-v),
\eeq
where we use the $R$-matrix \eqref{zamolr} and extend the subscript notation
$T^+_a(u)$ for $a=1,\dots,k$ from \eqref{tonettwo} to denote
the corresponding element of the
tensor product algebra
\ben
\underbrace{\End\CC^{M|2n}\ot\dots\ot\End\CC^{M|2n}}_k{}\ot\X^+(\g).
\een
Define the
{\em descending} filtration
on $\X^+(\g)$
by setting the degree of the generator $t_{ij}^{(-r)}$ with $r\geqslant 1$
to be equal to $r$. We let
$\wh\X^+(\g)$ denote the completion of $\X^+(\g)$
with respect to this filtration.

The {\em dual Yangian}
$\wh\Y^+(\g)$ is now defined as the quotient of the
algebra $\wh\X^+(\g)$ by the relations
$
T^{+}(u+\ka)^t\ts T^+(u)=1.
$
Consider the {\em ascending} filtration on $\wh\Y^+(\g)$
defined by $\deg t_{ij}^{(-r)}=-r$
for all $r\geqslant1$. A version of the Poincar\'e--Birkhoff--Witt theorem
holds for $\wh\Y^+(\g)$, which
can be proved by extending the approach of \cite{amr:rp}
to the super settings; see also \cite[Sec.~2.3]{k:op}.
It implies
the isomorphism for the associated graded algebra
\beql{grady}
\gr\wh\Y^+(\g)\cong \U\big(t^{-1}\g[t^{-1}]\big).
\eeq
The image $\bar t_{ij}^{\ts(-r)}$ of the generator
$t_{ij}^{(-r)}$ in the $(-r)$-th
component of the graded algebra $\gr\wh\Y^+(\g)$ corresponds to the element $F_{ij}[-r](-1)^{\bi}$
of $\U(t^{-1}\g[t^{-1}])$.

We will regard
$\wh\Y^+(\g)$ as a module over the Cartan subalgebra $\h$ of $\g$,
where each basis element $F_{ii}$ of $\h$ with $i=1,\dots,n+m$ acts as a derivation and
the action is defined on the generators by
\beql{fiiact}
F_{ii}\cdot t^+_{kl}(u)=\de_{ki}\ts t^+_{il}(u)- \de_{il}\ts t^+_{ki}(u)
-\de_{ki'}\ts t^+_{i'l}(u)+ \de_{il'}\ts t^+_{ki'}(u).
\eeq
Here we used the observation that the defining relations \eqref{RTT} and \eqref{RTTdualbcd}
written in terms of the generating series have the same form. The action \eqref{fiiact}
corresponds to the adjoint action of $F_{ii}$ on $t_{kl}(u)$ coming
from the embedding \eqref{embla}.
Denote by $\wh\Y^+(\g)^{\h}$ the
subalgebra of $\h$-invariants under this action.
Consider the left ideal $J$ of the algebra $\wh\Y^+(\g)$ generated by
all elements $t^{(-r)}_{ij}$ with the conditions
$1'\geqslant i>j\geqslant 1$ and $r\geqslant 1$.
The quotient of $\wh\Y^+(\g)^{\h}$ by the
two-sided ideal
$\wh\Y^+(\g)^{\h}\cap J$ is isomorphic
to the commutative superalgebra generated by the images
of the elements $t_{ii}^{(-r)}$ with $i=1,\dots,1'$ and $r\geqslant 1$
in the quotient. We will use the notation $\la^{(-r)}_i$
for this image.
An analogue of the Harish-Chandra
homomorphism \eqref{yhchbcd} now takes the form
\beql{yhchdualbcd}
\wh\Y^+(\g)^{\h}\to
\text{\quad\rm formal series in\quad}\la^{(-r)}_i,\quad i=1,\dots,1',\ r\geqslant 1,
\eeq
where the algebra of formal series is the completion of the algebra of polynomials
with respect to the gradation defined by setting the degree of $\la^{(-r)}_i$
to be equal to $r$.
We combine the elements $\la^{(-r)}_i$ into the formal series
\ben
\la^+_i(u)=1-\sum_{r=1}^{\infty} \la^{(-r)}_i u^{r-1},\qquad i=1,\dots,1',
\een
which will be understood as the images of the series
$t^+_{ii}(u)$ under the homomorphism \eqref{yhchdualbcd}.

To state the dual Yangian analogues of the results of Section~\ref{subsec:yc},
we will work with the tensor product algebra
\beql{tenprkydy}
\underbrace{\End\CC^{M|2n}\ot\dots\ot\End\CC^{M|2n}}_k{}\ot\wh\Y^+(\g).
\eeq
These analogues follow in a uniform way for odd and even values of $M$.
Introduce the formal series $\TT^{+\ts(k)}(u)$ with coefficients in
the dual Yangian $\wh\Y^+(\g)$ by the formula
\beql{bethebdy}
\TT^{+\ts(k)}(u)=\str_{1,\dots,k}\ts S^{(k)}T^+_1(u)\ts T^+_2(u+1)\dots T^+_k(u+k-1),
\eeq
with the supertrace taken over all $k$ copies of $\End\CC^{M|2n}$ in \eqref{tenprkydy}.

\bth\label{thm:hchbdy}
If $M=2m+1$, then
the image of the series $\TT^{+\ts(k)}(u)$ under the homomorphism \eqref{yhchdualbcd}
is given by
\ben
\TT^{+\ts(k)}(u)\mapsto\sum_{i_1>\dots>i_s\geqslant n'>i_{s+1}\geqslant\dots\geqslant i_p>n\geqslant i_{p+1}>\dots>i_k}
\la^+_{i_1}(u)\ts\la^+_{i_2}(u+1)\dots \la^+_{i_k}(u+k-1)(-1)^{\bi_1+\dots+\bi_k},
\een
assuming varying values of the parameters $s$ and $p$ in the summation,
with the condition that $n+m+1$ occurs among the summation indices $i_1,\dots,i_k$
at most once.
\eth

\bth\label{thm:hchdy}
If $M=2m$ and $m>n$, then
the image of the series $\TT^{+\ts(k)}(u)$ under the homomorphism \eqref{yhchdualbcd}
is given by
\ben
\TT^{+\ts(k)}(u)\mapsto\sum_{i_1>\dots>i_s\geqslant n'>i_{s+1}\geqslant\dots\geqslant i_p>n\geqslant i_{p+1}>\dots>i_k}
\la^+_{i_1}(u)\ts\la^+_{i_2}(u+1)\dots \la^+_{i_k}(u+k-1)(-1)^{\bi_1+\dots+\bi_k},
\een
assuming varying values of the parameters $s$ and $p$ in the summation,
with the condition that $n+m$ and $(n+m)'$ do not occur
simultaneously
among the summation indices $i_1,\dots,i_k$.
\eth

\bpf[Proof of Theorems~\ref{thm:hchbdy} and \ref{thm:hchdy}]
The theorems are immediate from the respective Corollaries~\ref{cor:hchbdy} and \ref{cor:hchdy}.
This follows since the defining
relations \eqref{RTTdualbcd} of the dual extended Yangian $\X^+(\g)$,
written in terms of the generating series $t^+_{ij}(u)$, take
the same form
as the Yangian relations \eqref{RTT} for the series $t_{ij}(u)$.
\epf

Our next step in proving Theorem~\ref{thm:hch} is to apply
Theorems~\ref{thm:hchbdy} and \ref{thm:hchdy} along with the isomorphism
\eqref{grady}. Recall that the vacuum module $V_{\cri}(\g)$ is isomorphic to
the universal enveloping algebra $\U(t^{-1}\g[t^{-1}])$ as a vector space.
We will produce elements of the dual Yangian whose graded images
will yield the Segal--Sugawara vectors $\Phi_k$; cf. \cite[Sec.~13.1]{m:so}.

Extend the ascending filtration on the dual Yangian
$\wh\Y^+(\g)$
defined by $\deg\ts t_{ij}^{(-r)}=-r$ to
the algebra of formal series $\wh\Y^+(\g)[[u,\di_u]]$ by setting
$\deg u=1$ and $\deg\di_u=-1$ so that the associated graded algebra is
isomorphic to $\U(t^{-1}\g[t^{-1}])[[u,\di_u]]$.
Consider the element of the extended algebra
\beql{onemtdbd}
\ga_k(M-2n)\ts\str^{}_{1,\dots,k}\tss S^{(k)} \big(1-T^+_1(u)\tss e^{-\di_u}\big)\dots
\big(1-T^+_k(u)\tss e^{-\di_u}\big),
\eeq
where we use the notation \eqref{ga}
and assume that $m>n$ in the case $M=2m$.
The element \eqref{onemtdbd}
has degree $-k$ and its image in the graded algebra coincides with
\beql{gradepbddiffo}
\ga_k(M-2n)\ts\str^{}_{1,\dots,k}\tss S^{(k)} \big(\di_u+\hF(u)_{+\ts 1}\big)\dots
\big(\di_u+\hF(u)_{+\ts k}\big),
\eeq
where the matrix $\hF(z)_+=[\hF_{ij}(z)_+]$ is defined in \eqref{fuplus}.

Now we will use Theorems~\ref{thm:hchbdy} and \ref{thm:hchdy} to evaluate
the image of the element \eqref{onemtdbd} under the natural extension
of the Harish-Chandra homomorphism \eqref{yhchdualbcd} to the algebra
$\wh\Y^+(\g)[[u,\di_u]]$, acting as the identity map on $u$ and $\di_u$.
Write \eqref{onemtdbd} in the form
\ben
\ga_k(M-2n)\ts\str^{}_{1,\dots,k}\ts S^{(k)}\ts \sum_{r=0}^k\sum_{1\leqslant i_1<\dots<i_r\leqslant k}
(-1)^r\ts T^+_{i_1}(u)\dots T^+_{i_r}(u-r+1)\tss e^{-r\tss\di_u}.
\een
Each product $T^+_{i_1}(u)\dots T^+_{i_r}(u-r+1)$
equals $P\ts T^+_{1}(u)\dots T^+_{r}(u-r+1)\ts P^{-1}$,
where $P$ is the image in \eqref{tenprkydy} (with the identity component
in the last tensor factor)
of a permutation $p\in\Sym_k$
such that $p(a)=i_a$ for $a=1,\dots,r$. Using the properties of the symmetrizer (see \cite[Sec.~1.2]{m:so})
and the cyclic property of supertrace, we can bring
the above expression to the form
\beql{epxa}
\ga_k(M-2n)\ts\sum_{r=0}^k (-1)^{r}\binom{k}{r} \ts
\str^{}_{1,\dots,k}\ts S^{(k)}\ts T^+_{1}(u)\dots T^+_{r}(u-r+1)\tss e^{-r\tss\di_u}.
\eeq
Now use \cite[Lemma~1.3.2]{m:so}
to calculate the partial supertraces of the symmetrizer $S^{(k)}$ over the copies
$r+1,\dots,k$ of $\End\CC^{M|2n}$ to get
\ben
\str^{}_{r+1,\dots,k}\ts S^{(k)}
=\frac{\ga_r(M-2n)}{\ga_k(M-2n)}\ts \binom{M-2n+k-2}{k-r}\ts\binom{k}{r}^{-1}\ts S^{(r)}.
\een
Then \eqref{epxa} equals
\ben
\sum_{r=0}^k (-1)^{r}\ts\ga_r(M-2n)\ts\binom{M-2n+k-2}{k-r} \ts
\str^{}_{1,\dots,r}\ts S^{(r)}\ts T^+_{1}(u)\dots T^+_{r}(u-r+1)\tss e^{-r\tss\di_u}.
\een
Applying the conjugation by the longest permutation in $\Sym_r$ and using the cyclic property
of supertrace we get
\ben
\str^{}_{1,\dots,r}\ts S^{(r)}\ts T^+_{1}(u)\dots T^+_{r}(u-r+1)
=\str^{}_{1,\dots,r}\ts T^+_{r}(u)\dots T^+_{1}(u-r+1)\ts S^{(r)}
\een
which equals
\ben
\str^{}_{1,\dots,r}\ts S^{(r)}\ts T^+_{1}(u-r+1)\dots T^+_{r}(u)
\een
by \eqref{symo} and the defining relations \eqref{RTTdualbcd} for the dual
Yangian.

The following property
of the noncommutative complete symmetric functions
\eqref{nonchl} will be used in the calculations below.

\ble\label{lem:comshi}
Let $u$ be
a complex parameter. We have the relation
\begin{multline}
h_k(-u+x_1,\dots,-u+x_n,u+x_{n+1},\dots,u+x_{(n+1)'},-u+x_{n'},\dots,-u+x_{1'})\\[0.4em]
{}=\sum_{r=0}^k\binom{M-2n+k-1}{k-r}\ts h_r(x_1,x_2,\dots,x_{1'})\ts u^{k-r}.
\label{hux}
\end{multline}
\ele

\bpf
Computing the generating function for the sequence on the right hand side, we get
\begin{multline}
\non
\sum_{k=0}^{\infty}\sum_{r=0}^k\binom{M-2n+k-1}{k-r}\tss h_r(x_1,\dots,x_{1'})\tss u^{k-r}\tss q^k\\
{}=\sum_{r=0}^{\infty}h_r(x_1,\dots,x_{1'})\tss q^r\sum_{k=r}^{\infty}
\binom{M-2n+k-1}{k-r}\tss u^{k-r}\tss q^{k-r}.
\end{multline}
This equals
\begin{multline}
\non
\sum_{r=0}^{\infty}h_r(x_1,\dots,x_{1'})\tss \frac{q^r}{(1-q\tss u)^{M-2n+r}}
=(1-q\tss u+q\tss x_1)\dots (1-q\tss u+q\tss x_n)\\[0.4em]
{}\times (1-q\tss u-q\tss x_{n+1})^{-1}\dots (1-q\tss u-q\tss x_{(n+1)'})^{-1}\ts
(1-q\tss u+q\tss x_{n'})\dots (1-q\tss u+q\tss x_{1'})
\end{multline}
and so coincides with the generating function of the sequence
on the left hand side of \eqref{hux}.
\epf

Now we will split the remaining arguments into two cases, depending on whether $M$ is odd or even.

\subsubsection*{Case $M=2m+1.$}

By Theorem~\ref{thm:hchbdy}, the Harish-Chandra image of the
expression \eqref{onemtdbd} is given by
\begin{multline}
\sum_{r=0}^k(-1)^{r}\ts\ga_r(M-2n)\ts\binom{M-2n+k-2}{k-r}\\[0.4em]
{}\times\sum_{i_1>\dots>i_s\geqslant n'>i_{s+1}\geqslant\dots\geqslant i_p>n\geqslant i_{p+1}>\dots>i_r}
\la^+_{i_1}(u-r+1)\dots \la^+_{i_r}(u)\tss e^{-r\tss\di_u}(-1)^{\bi_1+\dots+\bi_r}
\non
\end{multline}
with the condition that $n+m+1$ occurs among the summation indices $i_1,\dots,i_r$
at most once. We can rewrite the image in the form
\begin{multline}
\label{lahchb}
\sum_{r=0}^k(-1)^{r}\ts\ga_r(M-2n)\ts\binom{M-2n+k-2}{k-r}\\[0.4em]
{}\times\sum_{j_1<\dots<j_s\leqslant n<j_{s+1}\leqslant\dots\leqslant j_p<n'\leqslant j_{p+1}<\dots<j_r}
\la^+_{j_1}(u)\tss e^{-\di_u}\dots \la^+_{j_r}(u)\tss e^{-\di_u}(-1)^{\bj_1+\dots+\bj_r}
\end{multline}
with the condition that $n+m+1$ occurs among the summation indices $j_1,\dots,j_r$
at most once.

Recall the notation \eqref{muir} and combine the elements $\mu_i[-r]$ into the series
\beql{mui}
\mu_i(u)=\sum_{r=1}^{\infty}\mu_i[-r]\ts u^{r-1}.
\eeq
Set
\beql{sila}
\nu_i(u)=\big(1-\la^+_{i}(u)\tss e^{-\di_u}\big)(-1)^{\bi},\qquad i=1,\dots,1',
\eeq
and observe that each element $\nu_i(u)$ has degree $-1$,
with the top degree component equal to $(-1)^{\bi}\ts\di_u+\mu_i(u)$.
Our next goal is to express the image \eqref{lahchb} in terms
of the variables $\nu_i(u)$ and find its top degree component as the graded image
in the algebra $\U(t^{-1}\g[t^{-1}])[[u,\di_u]]$.

\ble\label{lem:lasib}
The expression \eqref{lahchb} multiplied by $2\ts(-1)^{k+1}\ts\binom{M/2-n-2}{M-2n+k-2}$
equals
\ben
\bal
\sum_{r=0}^k&(-1)^r\ts\binom{M/2-n-2}{M-2n+r-3}\\[0.4em]
{}&\times\ts
\sum_{a_1+\dots+a_{1'}=r}
\nu_{1}(u)^{a_1}\dots \nu_{n+m}(u)^{a_{n+m}}
\nu_{(n+m)'}(u)^{a_{(n+m)'}}\dots \nu_{1'}(u)^{a_{1'}}{}\\[1em]
{}+\sum_{r=1}^k&(-1)^r\ts\binom{M/2-n-2}{M-2n+r-3}\\[0.4em]
{}&\times\ts\sum_{a_1+\dots+a_{1'}=r-1}
\nu_{1}(u)^{a_1}\dots \nu_{n+m}(u)^{a_{n+m}}\big(\nu_{n+m+1}(u)-2\big)
\nu_{(n+m)'}(u)^{a_{(n+m)'}}\dots \nu_{1'}(u)^{a_{1'}},
\eal
\een
where $a_1,\dots,a_{1'}$ run over nonnegative integers and
each of $a_i,a_{i'}$
takes only two values $0$ and $1$ for $i=1,\dots,n$.
\ele

\bpf
The statement is verified by substituting \eqref{sila} into both terms
and calculating the coefficients of the sum
\beql{sumla}
\sum_{i_1<\dots<i_s\leqslant n<i_{s+1}\leqslant\dots\leqslant i_p<n'\leqslant i_{p+1}<\dots<i_r}
\la^+_{i_1}(u)\tss e^{-\di_u}\dots \la^+_{i_r}(u)\tss e^{-\di_u}(-1)^{\bi_1+\dots+\bi_r},
\eeq
with the condition that $n+m+1$ occurs among the summation indices $i_1,\dots,i_r$
at most once. To evaluate the coefficients, apply Lemma~\ref{lem:comshi}
with $u=1$ and $x_i=\la^+_{i}(u)\tss e^{-\di_u}(-1)^{\bi}$
to all expressions. Observe that the factor $\nu_{n+m+1}(u)-2$ becomes
$-1-x_{n+m+1}$, and that the relation of Lemma~\ref{lem:comshi}
and its proof extends to the counterparts of \eqref{nonchl} given by
\ben
\sum_{a_1+\dots+a_{1'}=l}
x^{a_1}_{1}\dots x^{a_{n+m}}_{n+m}\ts x_{n+m+1}\ts x^{a_{(n+m)'}}_{(n+m)'}\dots x^{a_{1'}}_{1'},
\een
where $a_1,\dots,a_{1'}$ run over nonnegative integers and each of $a_i,a_{i'}$
takes only two values $0$ and $1$ for $i=1,\dots,n$.
The coefficient of the sum \eqref{sumla}
in the resulting expression is given by the sum
\ben
\sum_{l=r}^k (-1)^{l-r}\binom{M/2-n-2}{M-2n+l-3}\binom{M-2n+l-3}{l-r}
\een
which is easily evaluated with the use of the simple identity
\ben
\sum_{k=0}^n(-1)^k\ts\binom{x}{k}=(-1)^n\ts\binom{x-1}{n},
\een
where $x$ is a variable. The sum equals
\ben
\binom{M/2-n-2}{M-2n+r-3}\binom{M/2-n+k-1}{k-r},
\een
which coincides with
\ben
2\ts (-1)^{k-r+1}\ts\ga_r(M-2n)\ts\binom{M/2-n-2}{M-2n+k-2}\ts\binom{M-2n+k-2}{k-r},
\een
as claimed.
\epf

Now the proof of Theorem~\ref{thm:hch} is completed by adapting the corresponding argument
of \cite[Sec.~13.1, Type~$B_n$]{m:so} to the super case.
Denote the expression in Lemma~\ref{lem:lasib} by $A_k$.
Since the degree of
the element \eqref{onemtdbd} is $-k$, its Harish-Chandra image \eqref{lahchb}
and the expression $A_k$ also have degree $-k$.
Observe that the terms in both sums of $A_k$ are independent of $k$ so that
$A_{k+1}=A_k+B_{k+1}$, where
\ben
\bal
B_{k+1}={}&(-1)^{k+1}\ts\binom{M/2-n-2}{M-2n+k-2}\\[0.4em]
{}&\times\ts
\sum_{a_1+\dots+a_{1'}=k+1}
\nu_{1}(u)^{a_1}\dots \nu_{n+m}(u)^{a_{n+m}}
\nu_{(n+m)'}(u)^{a_{(n+m)'}}\dots \nu_{1'}(u)^{a_{1'}}{}\\[1em]
{}+{}&(-1)^{k+1}\ts\binom{M/2-n-2}{M-2n+k-2}\\[0.4em]
{}&\times\ts\sum_{a_1+\dots+a_{1'}=k}
\nu_{1}(u)^{a_1}\dots \nu_{n+m}(u)^{a_{n+m}}\big(\nu_{n+m+1}(u)-2\big)
\nu_{(n+m)'}(u)^{a_{(n+m)'}}\dots \nu_{1'}(u)^{a_{1'}},
\eal
\een
where $a_1,\dots,a_{1'}$ run over nonnegative integers and
each of $a_i,a_{i'}$
takes only two values $0$ and $1$ for $i=1,\dots,n$.

Since $A_{k+1}$ has degree $-k-1$, its
component of degree $-k$ is zero, and so the sum of the homogeneous components
of degree $-k$ of $A_k$ and $B_{k+1}$ is zero.
However,
each element $\nu_i(u)$ has degree $-1$
and so
the
component of $A_k$ of degree $-k$ equals the component of degree $-k$
of the term
\ben
2\ts(-1)^{k+1}\ts\binom{M/2-n-2}{M-2n+k-2}\ts\sum_{a_1+\dots+a_{1'}=k}
\nu_{1}(u)^{a_1}\dots \nu_{n}(u)^{a_n}
\nu_{n'}(u)^{a_{n'}}\dots \nu_{1'}(u)^{a_{1'}}.
\een
Taking into account the constant factor used in Lemma~\ref{lem:lasib},
we can conclude that the Harish-Chandra image of \eqref{gradepbddiffo}
is given by the
noncommutative complete symmetric function
\begin{multline}\label{hchser}
h_k\big({-}\di_u+\mu_{1}(u),\dots,-\di_u+\mu_{n}(u),\di_u+\mu_{n+1}(u),\dots,\di_u+\mu_{n+m}(u),\\[0.4em]
\di_u+\mu_{(n+m)'}(u),\dots,\di_u+\mu_{(n+1)'}(u),
-\di_u+\mu_{n'}(u),\dots,-\di_u+\mu_{1'}(u)\big).
\end{multline}
The proof of Theorem~\ref{thm:hch} in the case $M=2m+1$
is completed by evaluating the coefficients of
the polynomials in $\di_u$ which appear in
\eqref{gradepbddiffo} and \eqref{hchser} at $u=0$. The constant terms
of these polynomials are the Segal--Sugawara vector $\Phi_k$ and its
Harish-Chandra image $\f(\Phi_k)$, respectively.
The evaluation at $u=0$ relies on the vertex algebra structure on
the vacuum module $V_{\cri}(\g)\cong\U(t^{-1}\g\tss[t^{-1}])$
and is equivalent to the application of the vacuum axiom, since
the constant term in \eqref{gradepbddiffo} is obtained by
the application of the Sugawara operator \eqref{phiksugasym}
to the vacuum vector.

\subsubsection*{Case $M=2m$ with $m>n$.}

The assumption $m>n$ ensures that the symmetrizer $S^{(k)}$
given by \eqref{symo}
is well-defined for all $k$. Therefore, to prove
Theorem~\ref{thm:hch} under this assumption,
we can use the same argument as in the case $M=2m+1$.
By Theorem~\ref{thm:hchdy}, the Harish-Chandra image of the
expression \eqref{onemtdbd} is given by
\begin{multline}
\label{lahchd}
\sum_{r=0}^k(-1)^{r}\ts\ga_r(2m-2n)\ts\binom{2m-2n+k-2}{k-r}\\[0.4em]
{}\times\sum_{j_1<\dots<j_s\leqslant n<j_{s+1}\leqslant\dots\leqslant j_p<n'\leqslant j_{p+1}<\dots<j_r}
\la^+_{j_1}(u)\tss e^{-\di_u}\dots \la^+_{j_r}(u)\tss e^{-\di_u}(-1)^{\bj_1+\dots+\bj_r}
\end{multline}
with
the condition that $n+m$ and $(n+m)'$ do not occur
simultaneously
among the summation indices $i_1,\dots,i_k$.
Introducing new variables by the same formulas \eqref{sila}
we come to the following counterpart of Lemma~\ref{lem:lasib},
where we use the notation
\ben
c_k=-\binom{2\tss m-2\tss n+k-2}{m-n-1}^{-1}.
\een
It takes essentially the same form as \cite[Lemma~13.1.11]{m:so} for the case
of even orthogonal Lie algebras.

\ble\label{lem:lasid}
The expression \eqref{lahchd} multiplied by $2\tss c_k$ equals
\begin{multline}
2\tss c_k\ts
\sum_{\underset{\scriptstyle a^{}_{n+m}=a_{(n+m)'}=0}{a_1+\dots+a_{1'}=k}}
\nu_{1}(u)^{a_1}\dots \nu_{1'}(u)^{a_{1'}}\\
{}+c_k\ts
\sum_{\underset{\scriptstyle \text{only one of\ } a^{}_{n+m}\text{\ and\ }a_{(n+m)'}
\text{\ is zero}}{a_1+\dots+a_{1'}=k}}
\nu_{1}(u)^{a_1}\dots \nu_{1'}(u)^{a_{1'}}\\[0.5em]
{}-\sum_{r=1}^k\ts\frac{r\ts c_r}{m-n+r-1}\ts
\sum_{\underset{\scriptstyle a^{}_{n+m}=a_{(n+m)'}=0}{a_1+\dots+a_{1'}=r}}
\nu_{1}(u)^{a_1}\dots \nu_{1'}(u)^{a_{1'}}\\[0.5em]
{}+\sum_{r=1}^k\ts\frac{(m-n-1)\ts c_r}{m-n+r-1}\ts
\sum_{\underset{\scriptstyle \text{only one of\ } a^{}_{n+m}\text{\ and\ }a_{(n+m)'}
\text{\ is zero}}{a_1+\dots+a_{1'}=r}}
\nu_{1}(u)^{a_1}\dots \nu_{1'}(u)^{a_{1'}},
\non
\end{multline}
where $a_1,\dots,a_{1'}$ run over nonnegative integers and
each of $a_i,a_{i'}$
takes only two values $0$ and $1$ for $i=1,\dots,n$.
\ele

\bpf
Set
$x_i=\la^+_{i}(u)\tss e^{-\di_u}(-1)^{\bi}$ so that
\eqref{sila} takes the form
\beql{silad}
\nu_i(u)=(-1)^{\bi}-x_i,\qquad i=1,\dots,1'.
\eeq
Substitute \eqref{silad} into the expression given in the lemma
and calculate the coefficients of the monomials $x_{j_1}\dots x_{j_r}$
with the conditions on the subscripts as in \eqref{lahchd}.
To perform the calculation,
consider three families of noncommutative complete symmetric functions
$h_l(\nu^0),h_l(\nu^+)$ and $h_l(\nu^-)$,
defined in \eqref{nonchl}, depending on the respective
sets of variables
\ben
\bal
\nu^0&=\big(\nu_1(u),\dots,\nu_{n+m-1}(u),\nu_{(n+m-1)'}(u),\dots,\nu_{1'}(u)\big),\\
\nu^+&=\big(\nu_1(u),\dots,\nu_{n+m-1}(u),\nu_{n+m}(u),\nu_{(n+m-1)'}(u),\dots,\nu_{1'}(u)\big),\\
\nu^-&=\big(\nu_1(u),\dots,\nu_{n+m-1}(u),\nu_{(n+m)'}(u),\nu_{(n+m-1)'}(u),\dots,\nu_{1'}(u)\big).
\eal
\een
Using this notation, we can write the expression in the lemma in the form
\beql{compfo}
\bal
c_k\ts \big(h_k(\nu^+)+h_k(\nu^-)\big)
{}&+\sum_{r=1}^k\ts\frac{(m-n-1)\ts c_r}{m-n+r-1}\ts \big(h_r(\nu^+)+h_r(\nu^-)\big)\\
{}&-\sum_{r=1}^k\ts\frac{(2m-2n+r-2)\ts c_r}{m-n+r-1}\ts h_r(\nu^0).
\eal
\eeq
Now apply Lemma~\ref{lem:comshi} with $u=1$ to express $h_l(\nu^0),h_l(\nu^+)$ and $h_l(\nu^-)$
as the respective linear combinations of the noncommutative complete symmetric functions
$h_r(x^0),h_r(x^+)$ and $h_r(x^-)$, where
\ben
\bal
x^0&=\big(x_1,\dots,x_{n+m-1},x_{(n+m-1)'},\dots,x_{1'}\big),\\
x^+&=\big(x_1,\dots,x_{n+m-1},x_{n+m},x_{(n+m-1)'},\dots,x_{1'}\big),\\
x^-&=\big(x_1,\dots,x_{n+m-1},x_{(n+m)'},x_{(n+m-1)'},\dots,x_{1'}\big).
\eal
\een
Then \eqref{compfo} takes the form
\ben
\bal
c_k\ts\sum_{s=0}^k(-1)^s&\binom{2m-2n+k-2}{k-s}\big(\wt h_s(x)+h_s(x^0)\big)\\[0.4em]
{}+&\sum_{r=1}^k\ts\frac{(m-n-1)\ts c_r}{m-n+r-1}\ts
\sum_{s=0}^r(-1)^s\binom{2m-2n+r-2}{r-s}\big(\wt h_s(x)+h_s(x^0)\big)\\[0.4em]
{}-&\sum_{r=1}^k\ts\frac{(2m-2n+r-2)\ts c_r}{m-n+r-1}\ts
\sum_{s=0}^r(-1)^s\binom{2m-2n+r-3}{r-s}\ts h_s(x^0),
\eal
\een
where we set $\wt h_s(x)=h_s(x^+)+h_s(x^+)-h_s(x^0)$. For a given $s\in\{1,\dots,k\}$,
the coefficient of $(-1)^s\tss h_s(x^0)$ in this expression equals
\ben
\bal
c_k\ts\binom{2m-2n+k-2}{k-s}
{}&+\sum_{r=s}^k\ts\frac{(m-n-1)\ts c_r}{m-n+r-1}\ts \binom{2m-2n+r-2}{r-s}\\[0.4em]
{}&-\sum_{r=s}^k\ts\frac{(2m-2n+r-2)\ts c_r}{m-n+r-1}\ts
\binom{2m-2n+r-3}{r-s}.
\eal
\een
This coefficient is zero; this follows easily as an application of the identity
\ben
\sum_{i=0}^k\ts\binom{y+i}{i}=\binom{y+k+1}{k},
\een
where $y$ is a variable. Essentially the same calculation
of the coefficient of $(-1)^s\tss\wt h_s(x)$
leads to the conclusion
that \eqref{compfo} equals
\ben
2\tss c_k\ts\sum_{r=0}^k(-1)^{r}\ts\ga_r(M-2n)\ts\binom{M-2n+k-2}{k-r}\ts \wt h_r(x),
\een
which coincides with \eqref{lahchd} multiplied by $2\tss c_k$,
as required.
\epf

Let $A_k$ denote the four-term expression in Lemma~\ref{lem:lasid}. This expression
equals $2\tss c_k$ times the Harish-Chandra image of \eqref{onemtdbd} and so
$A_k$ has degree $-k$. Hence, the component of degree $-k$
of the expression $A_{k+1}$ is zero.
On the other hand,
each element $\nu_i(u)$ has degree $-1$
with the top degree component equal to $(-1)^{\bi}\di_u+\mu_i(u)$,
where $\mu_i(u)$ is defined in \eqref{mui}.
This implies that the component of degree $-k$ in the sum
of the third and fourth terms in $A_k$ is zero. Therefore, the
component of $A_k$ of degree $-k$ equals the component of degree $-k$
in the sum of the first and the second terms.
Taking into account the constant factor $2\tss c_k$,
we conclude that the Harish-Chandra image of \eqref{gradepbddiffo}
equals the component of degree $-k$
of the sum
\ben
\sum_{\underset{\scriptstyle a^{}_{n+m}=a_{(n+m)'}=0}{a_1+\dots+a_{1'}=k}}
\nu_{1}(u)^{a_1}\dots \nu_{1'}(u)^{a_{1'}}
+\frac12\ts
\sum_{\underset{\scriptstyle \text{only one of\ } a^{}_{n+m}\text{\ and\ }a_{(n+m)'}
\text{\ is zero}}{a_1+\dots+a_{1'}=k}}
\nu_{1}(u)^{a_1}\dots \nu_{1'}(u)^{a_{1'}},
\een
and hence coincides with the half-sum of two
noncommutative complete symmetric functions
\begin{multline}\non
h_k\big({-}\di_u+\mu_{1}(u),\dots,-\di_u+\mu_{n}(u),\di_u+\mu_{n+1}(u),\dots,\di_u+\mu_{n+m}(u),\\[0.4em]
\di_u+\mu_{(n+m-1)'}(u),\dots,\di_u+\mu_{(n+1)'}(u),
-\di_u+\mu_{n'}(u),\dots,-\di_u+\mu_{1'}(u)\big)
\end{multline}
and
\begin{multline}\non
h_k\big({-}\di_u+\mu_{1}(u),\dots,-\di_u+\mu_{n}(u),\di_u+\mu_{n+1}(u),\dots,\di_u+\mu_{n+m-1}(u),\\[0.4em]
\di_u+\mu_{(n+m)'}(u),\dots,\di_u+\mu_{(n+1)'}(u),
-\di_u+\mu_{n'}(u),\dots,-\di_u+\mu_{1'}(u)\big).
\end{multline}
As in the case of odd $M$,
the proof of Theorem~\ref{thm:hch} for $M=2m$ with $m>n$
is completed by evaluating the coefficients of
the polynomials in $\di_u$ that appear in
\eqref{gradepbddiffo} and in this half-sum at $u=0$.

\subsection{Extrapolation argument}
\label{subsec:ea}

Here we complete the proof of Theorem~\ref{thm:hch} in the remaining case $M=2m$ with
the conditions $1\leqslant m\leqslant n$. Although the Segal--Sugawara vectors $\Phi^{}_{k}$ given in
\eqref{phimliesup} are defined for all values of the parameters, the argument
in the previous section relies on the use of the symmetrizer $S^{(k)}$
and is only valid in the region $m> n$. To extrapolate the formulas to the region
$1\leqslant m\leqslant n$,
we will fix values of $n$ and $k$, and let $m$ vary, assuming the natural
embeddings of the orthosymplectic Lie superalgebras
\beql{chainosp}
\osp_{0|2n}\subset \osp_{2|2n}\subset \dots\subset \osp_{2m|2n}\subset \osp_{2m+2|2n}\subset\dots,
\eeq
where $F_{ij}\in \osp_{2m|2n}$ is identified with the element of $\osp_{2m+2|2n}$
with the same name.
We will analyse the Harish-Chandra
images $\f(\Phi^{}_{k})$, regarding the coefficients of the monomials
in the variables $\mu_i[r]$ as functions in $m$. Our goal is to demonstrate
that these coefficients are polynomials in $m$ of certain fixed degrees.
Since we know these polynomials for the infinite set of values $m>n$, we will be able to conclude
that the formulas for the images $\f(\Phi^{}_{k})$ are valid for all values of $m\geqslant 1$.

First, observe that the coefficient $\Yc_{k,\ell}(2m-2n-1)$ occurring in \eqref{phimliesup}
is a polynomial in $m$ of degree $k-l$. Therefore, it will be enough to
establish the desired polynomiality property for the Harish-Chandra
images of the elements $\str^{}_{1,\dots,\ell}\ts H^{(\ell)} \hF[-\la]$.
In their turn, these elements are linear combinations of the supertraces
$\str^{}_{1,\dots,\ell}\ts P_{\si} \hF[-\la]$, where $P_{\si}$ is the image
of a permutation $\si\in\Sym_l$ under the action of $\Sym_l$ on the
superspace $(\CC^{2m|2n})^{\ot l}$.

To make a further reduction of the family of elements
whose Harish-Chandra
images would be sufficient to consider, use the Brauer algebra action
\eqref{brahom}, where we take $M=2m$ and $\om=2m-2n$.
As we will work with the extended tensor product superalgebra,
we will usually identify these elements with
$P_{ab} \otimes 1$ and $Q_{ab} \otimes 1$ in \eqref{tenprka}, respectively.

For any diagram $d\in \Bc_k(\om)$ we let $D$ denote its image under the homomorphism \eqref{brahom}.
Given any negative integers $a_1,\dots,a_k$, introduce elements of \eqref{tenprka} by
\ben
\FF_{ij}(a_1,\dots,a_k)=\big(\hF[a_1]\dots \hF[a_k]\big)_{ij}
\een
and set
\beql{xadedef}
X(a_1,\dots,a_k)=\str\ts \hF[a_1]\dots \hF[a_k]=\sum_{i=1}^{1'} \ts \FF_{ii}(a_1,\dots,a_k)(-1)^{\bi}.
\eeq
We assume that $\FF_{ij}(a_1,\dots,a_k)=\de_{ij}$ for $k=0$.

\ble\label{lem:redu}
The element
\beql{trdf}
\str^{}_{1,\dots,k} D \hF[a_1]_1\dots \hF[a_k]_k
\eeq
of the superalgebra \eqref{tenprka}
can be written as a linear combination of products of the form
\beql{proxx}
X(b_1,\dots,b_r)\dots X(c_1,\dots,c_s),
\eeq
whose coefficients are polynomials in $m$ of degrees bounded by $k$, where
each of $b_1,\dots,c_s$ is a sum of some numbers $a_i$.
\ele

\bpf
Recall that $t_l$ denotes the partial super-transposition
\eqref{suptra} on the superalgebra \eqref{tenprka}, as well as the partial
transposition on the Brauer algebra $\Bc_k(\om)$. This coincidence is unambiguous
because the operations are consistent in the sense that the image
of a transposed diagram $d^{\tss t_l}$ under the homomorphism \eqref{brahom}
coincides with the operator $D^{t_l}$.

As the first step of the proof, observe that for any $l\in\{1,\dots,k\}$ we have the identity
\beql{srtdddta}
\str^{}_{1,\dots,k} (D+D^{t_l}) \hF[a_1]_1\dots \hF[a_k]_k=0.
\eeq
This is clear from the property
\ben
\str^{}_{1,\dots,k} XY= \str^{}_{1,\dots,k} X^{t_l}Y^{t_l}
\een
which holds for arbitrary elements $X$ and $Y$ of \eqref{tenprka}.
The identity \eqref{srtdddta} follows by taking $X=D+D^{t_l}$ and $Y=\hF[a_1]_1\dots \hF[a_k]_k$
and noting that $\hF[a_l]^t=-\hF[a_l]$ by \eqref{symma}.

It is an easy consequence of \eqref{srtdddta} that, up to a sign,
the element \eqref{trdf}
equals an element of the same form, where $D$ is a diagram without horizontal edges.
Hence, we may assume that $D=P_{\si}$ for a permutation $\si\in\Sym_k$.

Furthermore, if $i_1,\dots,i_r$ are distinct elements of the set $\{1,\dots,k\}$
and $\si=(i_r,\dots,i_1)$ is the associated cycle, then
\beql{srtpsx}
\str^{}_{i_1,\dots,i_r} P_{\si} \hF[b_1]_{i_1}\dots \hF[b_r]_{i_r}=X(b_1,\dots,b_r)
\eeq
for any negative integers $b_1,\dots,b_r$. Indeed,
we can write $P_{\si}=P_{i_{r-1}i_r}\dots P_{i_1i_2}$ so that
\ben
\bal
\str^{}_{i_1,\dots,i_r} P_{\si} \hF[b_1]_{i_1}\dots \hF[b_r]_{i_r}
{}&=\str^{}_{i_1,\dots,i_r} P_{i_{r-1}i_r}\dots P_{i_2i_3}\hF[b_1]_{i_2}P_{i_1i_2} \dots \hF[b_r]_{i_r}\\[0.4em]
{}&=\str^{}_{i_2,\dots,i_r} P_{i_{r-1}i_r}\dots P_{i_2i_3}\hF[b_1]_{i_2}\hF[b_2]_{i_2} \dots \hF[b_r]_{i_r}
\eal
\een
and \eqref{srtpsx} follows by an easy induction, where we used the relation $\str_a\ts P_{ab}=1$.

To complete the proof of the lemma, write an arbitrary element $\si\in\Sym_k$
as a product of disjoint cycles. To apply \eqref{srtpsx} to the supertrace
\beql{pstrdf}
\str^{}_{1,\dots,k} P_{\si} \hF[a_1]_1\dots \hF[a_k]_k,
\eeq
we will reorder the factors $\hF[a_i]_i$ to match the cycle decomposition of $\si$
by using the commutation relations
\beql{comhatf}
\big[\hF[a]_i,\hF[b]_j\big]=(P_{ij}-Q_{ij})\hF[a+b]_j-\hF[a+b]_j(P_{ij}-Q_{ij}),
\eeq
implied by \eqref{comrellie}. Arguing now by induction on $k$, we get
a formula for \eqref{pstrdf} as a product \eqref{proxx} plus
a linear combination of
expressions of the form \eqref{trdf} with smaller values of $k$.
The proof is completed by using its first part and the induction hypothesis.
The coefficients of the linear combination are polynomials in $\om=2m-2n$
arising from the partial supertraces $\str^{}_l\ts D$ of the operators representing Brauer
diagrams. It is clear from the above argument that the degree of $m$ in such polynomials
does not exceed $k$.
\epf

By Lemma~\ref{lem:redu}, it is now sufficient to establish the polynomiality property
for the Harish-Chandra images of the products of the form \eqref{proxx}.
Since each factor lies in the centralizer $\U(t^{-1}\g[t^{-1}])^{\h}$,
it is enough to consider the Harish-Chandra image of each factor
due to the homomorphism property of $\f$ in \eqref{hchaff}.
Moreover, taking into account \eqref{xadedef}, we may reduce this further to checking
the polynomiality property for the Harish-Chandra images of the elements
$\FF_{ii}(a_1,\dots,a_k)$ for all $i=1,\dots,1'$.

\ble\label{lem:cof}
We have the commutation relations
\ben
\bal
\big[\hF_{ij}[a],\FF_{kl}(b_1,\dots,b_r)\big]
{}&=\sum_{s=1}^r
\big(\FF_{kj}(b_1,\dots,b_{s-1})\ts\FF_{il}(a+b_s,\dots,b_{r})\\
{}&\qquad-\FF_{kj}(b_1,\dots,a+b_{s})\ts\FF_{il}(b_{s+1},\dots,b_{r})\big)(-1)^{\bi\bj+\bi\bk+\bj\bk}\\
{}&-\sum_{s=1}^r
\big(\FF_{ki'}(b_1,\dots,b_{s-1})\ts\FF_{j'l}(a+b_s,\dots,b_{r})\\
{}&\qquad-\FF_{ki'}(b_1,\dots,a+b_{s})\ts\FF_{j'l}(b_{s+1},\dots,b_{r})\big)
(-1)^{\bj+\bi\bk+\bj\bk}\ta_i\ta_j,
\eal
\een
where $a,b_1,\dots,b_r$ are negative integers.
\ele

\bpf
Note that the particular case $r=1$ coincides
with the relation
\eqref{comhatf} written in terms
of matrix elements:
\ben
\bal
\big[\hF_{ij}[a],\hF_{kl}[b]\big]&=\big(\de_{kj}\hF_{il}[a+b]
- \de_{il}\hF_{kj}[a+b]\big)(-1)^{\bi\bj+\bi\bk+\bj\bk}\\[0.4em]
{}&- \big(\de_{ki'} \hF_{j'l}[a+b]
-\de_{j'l} \hF_{ki'}[a+b]\big)(-1)^{\bj+\bi\bk+\bj\bk}\ta_i\ta_j.
\eal
\een
The formula for the general values of $r$ then easily follows.
\epf

Now suppose that $i\leqslant n+m$ and consider the elements $\FF_{ii}(a_1,\dots,a_k)$.
To evaluate their Harish-Chandra images, write
\ben
\FF_{ii}(a_1,\dots,a_k)=\sum_{j=1}^{1'}\hF_{ij}[a_1]\ts \FF_{ji}(a_2,\dots,a_k).
\een
Hence,
\ben
\f\big(\FF_{ii}(a_1,\dots,a_k)\big)=\mu_i[a_1]\ts \f\big(\FF_{ii}(a_2,\dots,a_k)\big)(-1)^{\bi}
+\f\tss\big(\sum_{j<i}\hF_{ij}[a_1]\ts \FF_{ji}(a_2,\dots,a_k)\big).
\een
Now apply Lemma~\ref{lem:cof} to commute $\hF_{ij}[a_1]$ with $\FF_{ji}(a_2,\dots,a_k)$.
Since $j'>i$, the second sum in the commutator formula vanishes
after the application of the Harish-Chandra homomorphism and we get the recurrence formula
\ben
\bal
\f\big(\FF_{ii}(a_1,\dots,a_k)\big)&=\mu_i[a_1]\ts \f\big(\FF_{ii}(a_2,\dots,a_k)\big)(-1)^{\bi}\\[0.4em]
{}&+\sum_{j=1}^{i-1}\sum_{s=2}^k
\Bigg(\f\big(\FF_{jj}(a_2,\dots,a_{s-1})\big)\f\big(\FF_{ii}(a_1+a_s,\dots,a_k)\big)\\
{}&\qquad\qquad-\f\big(\FF_{jj}(a_2,\dots,a_1+a_{s})\big)
\f\big(\FF_{ii}(a_{s+1},\dots,a_k)\big)\Bigg)(-1)^{\bj}.
\eal
\een

To get a similar formula for the Harish-Chandra images
of the diagonal matrix elements $\FF_{i'i'}(a_1,\dots,a_k)$ with
$i\leqslant n+m$, note that these elements coincide with the $(i,i)$-entries
of the transposed matrix $(\hF[a_1]\dots \hF[a_k])^t$. Therefore,
the desired recurrence formula is implied by the following lemma.

\ble\label{lem:trans}
The transposed matrix $(\hF[a_1]\dots \hF[a_k])^t$ equals a linear combination
of expressions of the form
\ben
X(b_1,\dots,b_r)\dots  X(c_1,\dots,c_s) \hF[d_1]\dots \hF[d_l]
\een
with $l\leqslant k$,
where each of the symbols $b_1,\dots,d_l$ is a sum $a_{i_1}+\dots+a_{i_v}$, while the coefficients are
polynomials in $\om=2m-2n$ of degrees not exceeding $k$.
\ele

\bpf
We use induction on $k$ and work in the superalgebra \eqref{tenprka}
with two copies of $\End \CC^{2m|2n}$. We have
\ben
\bal
\big(\hF[a_1]_2\dots \hF[a_k]_2\big)^{t_2}&=\str_1 Q_{12}\big(\hF[a_1]_2\dots \hF[a_k]_2\big)^{t_2}
=\str_1 Q_{12}\hF[a_1]_1\dots \hF[a_k]_1\\
&=\str_1 Q_{12}\big(\hF[a_1]_2\dots \hF[a_{k-1}]_2\big)^{t_2} \hF[a_k]_1,
\eal
\een
where we used the relations $Q_{12}Y^{\tss t}_2=Q_{12}Y_1$ valid for arbitrary elements
$Y$ of \eqref{tenprka}, and the property $\str_1 Q_{12}=1$.
By the induction hypothesis, $\big(\hF[a_1]_2\dots \hF[a_{k-1}]_2\big)^{t_2}$ is a linear
combination of the terms listed in the lemma. Take one of these terms
and denote by $X$ the product of supertraces of the form
$X(b_1,\dots,b_r)$ occurring as factors in this term.
Now evaluate the supertrace
\ben
\str_1 Q_{12}X\ts \hF[d_1]_2\dots \hF[d_l]_2 \hF[a_k]_1=\str_1 X Q_{12}\ts \hF[d_1]_2\dots \hF[d_l]_2 \hF[a_k]_1.
\een
Use the commutation relations to move $\hF[a_k]_1$ to the left. We have
\beql{cr}
\big[\hF[d_i]_2,\hF[a_k]_1\big]=\hF[a_k+d_i]_2(P_{12}-Q_{12})-(P_{12}-Q_{12})\hF[a_k+d_i]_2.
\eeq
Since $Q_{12}\hF[a_k]_1=-Q_{12}\hF[a_k]_2$, it remains to evaluate
the contribution of the terms arising from the commutation; they are
of the form
\beql{paq}
\str_1 X Q_{12}\ts \hF[d_1]_2\dots \hF[d_{i-1}]_2(P_{12}-Q_{12})\hF[a_k+d_i]_2\dots \hF[d_l]_2
\eeq
with $i=1,\dots,l$;
the first expression on the right hand side of \eqref{cr} contributes a similar term.
Since $Q_{12}P_{12}=Q_{12}$ and $Q_{12}Y_2Q_{12}=(\str\ts Y)\ts Q_{12}$
for an arbitrary element $Y\in\End \CC^{2m|2n} \otimes \U$,
expression \eqref{paq}
equals the sum
\begin{multline}
\str_1 X Q_{12}\ts \hF[d_1]_1\dots \hF[d_{i-1}]_1\hF[a_r+d_i]_2\dots \hF[d_l]_2\\[0.4em]
-\str_1 X\cdot X(d_1,\dots,d_{i-1}) Q_{12}\hF[a_k+d_i]_2\dots \hF[d_l]_2.
\non
\end{multline}
The first term equals
\ben
\str_1 X Q_{12}\ts \big(\hF[d_1]_2\dots \hF[d_{i-1}]_2\big)^{t_2}\hF[a_k+d_i]_2\dots \hF[d_l]_2
\een
so that the induction hypothesis applies to get the resulting linear combination.
The second term already has the required form since $\str_1 Q_{12}=1$.
As the calculations show, the numerical coefficients of the resulting
linear combination representing the transposed matrix $(\hF[a_1]\dots \hF[a_k])^t$
can only arise from the relations $Q_{12}^2=\om\ts Q_{12}$
and $\str\ts 1=\om$; the maximal power of $\om$ will not exceed $k$.
\epf

Returning to the Segal--Sugawara vectors, we will use
a superscript to indicate the dependence on $m$ by writing $\Phi^{(m)}_k$
for the vector $\Phi_k$ associated with $\osp_{2m|2n}$. Fix $1\leqslant m_0\leqslant n$
and apply the inductive procedure of calculating the
Harish-Chandra image $\f(\Phi^{(m_0)}_k)$, as described by Lemmas~\ref{lem:redu}--\ref{lem:trans}.
Write this image
as a linear combination of the basis monomials
\beql{basmon}
\mu_{j_1}[r_1]\dots \mu_{j_s}[r_s],
\eeq
where $1\leqslant j_1\leqslant\dots\leqslant j_s\leqslant n+m_0$ and $r_1,\dots,r_s$
are negative integers.

Using the chain of embeddings \eqref{chainosp},
let $m$ vary taking all values $m\geqslant m_0$ and write the
Harish-Chandra image $\f(\Phi^{(m)}_k)$ in a similar way, as a linear combination of the
corresponding basis monomials.
The above arguments imply that
the coefficient of any monomial \eqref{basmon} in this linear combination
is a polynomial in $m$ whose degree does not exceed $k$.

On the other hand, the claimed image $\f(\Phi^{(m)}_k)$ in Theorem~\ref{thm:hch}
also has the polynomiality property. Indeed, for any $m\geqslant m_0$, the coefficients
of the monomials of the form \eqref{basmon} in that image
are found by setting $\mu_i[r]=0$ for all $n+m_0<i\leqslant n+m$ and $r<0$.
It is easily seen that
as a result of this evaluation we get the expression
\begin{multline}\label{evalmuz}
\sum_{a_1+\dots+a_{1'}+b=k}\binom{2m-2m_0+b-2}{b}\ts
\big({-}\tau+\mu_{1}[-1]\big)^{a_1}\dots \big(\tau+\mu_{n+m_0}[-1]\big)^{a_{n+m_0}}\tss \tau^b\\[0.4em]
{}\times \big(\tau+\mu_{(n+m_0)'}[-1]\big)^{a_{(n+m_0)'}}\dots \big({-}\tau+\mu_{1'}[-1]\big)^{a_{1'}}\ts 1,
\end{multline}
where $a_1,\dots,a_{1'}$ and $b$ run over nonnegative integers and each of $a_i,a_{i'}$
takes only two values $0$ and $1$ for $i=1,\dots,n$.
Hence the dependence on $m$ in the coefficient in question comes only from the binomial coefficients
and so it is a polynomial in $m$ whose degree does not exceed $k$.

By the arguments of Sec.~\ref{subsec:dy},
the coefficients
of the monomials \eqref{basmon}
in the expansion of $\f(\Phi^{(m)}_k)$
are known for all $m>n$ and given by the formulas of Theorem~\ref{thm:hch}.
Therefore, by the polynomiality property, the
coefficients are given by the same formulas for all values of $m\geqslant m_0$,
thus completing the proof of Theorem~\ref{thm:hch} in the remaining case $M=2m$
for all $m\geqslant 1$.

\subsection{New proof for symplectic Lie algebras}
\label{subsec:np}

In the above proof of Theorem~\ref{thm:hch}, the assumption $M\geqslant 1$ was used.
The case $M=0$ corresponds to the symplectic Lie algebras $\osp_{0|2n}\cong\spa_{2n}$,
and the Harish-Chandra images of the Segal--Sugawara vectors
$\Phi^{}_{k}$ are already known in this case; see \cite[Sec.~13.1]{m:so}.
The formula \eqref{phimliesup} was found in \cite{m:ss}, while its new proof
is given in \cite{mn:ss}. It is clear from the formula that
$\Phi_k=0$ for $k\geqslant 2n+2$.
We will show that the arguments of Sec.~\ref{subsec:ea}
extend to the symplectic case, thus producing a new proof of the formula for
the Harish-Chandra image of $\Phi^{}_{k}$.

Recall a noncommutative version of the elementary symmetric functions
in the variables $x_1,\dots,x_N$
defined by
\beql{elem}
e_k(x_1,\dots,x_N)=\sum_{i_1>\dots> i_k}
x_{i_1}\dots x_{i_k},
\eeq
for $k\geqslant 1$.

\bth\label{thm:sphch}
For any $1\leqslant k\leqslant 2n+1$,
the Harish-Chandra image $\f(\Phi_k)$ equals
\ben
(-1)^k\ts e_k\big(\tau+\mu_{1}[-1],\dots,
\tau+\mu_{n}[-1],\tau,\tau-\mu_{n}[-1],\dots
\tau-\mu_{1}[-1]\big)\ts 1.
\een
\eth

\bpf
Take $m_0=0$ in the arguments of Sec.~\ref{subsec:ea}, and note that by \eqref{evalmuz},
the evaluation of the Harish-Chandra image $\f(\Phi^{(m)}_k)$ at
$\mu_i[r]=0$ for all $n<i\leqslant n+m$ and $r<0$ yields the expression
\begin{multline}\label{evalmuzer}
\sum_{a_1+\dots+a_{1'}+b=k}\binom{2m+b-2}{b}\ts
\big({-}\tau+\mu_{1}[-1]\big)^{a_1}\dots \big({-}\tau+\mu_{n}[-1]\big)^{a_n}\tss \tau^b\\[0.4em]
{}\times \big({-}\tau+\mu_{n'}[-1]\big)^{a_{n'}}\dots \big({-}\tau+\mu_{1'}[-1]\big)^{a_{1'}}\ts 1.
\end{multline}
By the polynomiality property of the Harish-Chandra image $\f(\Phi^{(m)}_k)$
established in Sec.~\ref{subsec:ea},
it suffices to take $m=0$ in this expression. As a polynomial in $m$,
the binomial coefficient equals
\ben
\frac{(2m-1)(2m)\dots (2m+b-2)}{b\tss !}
\een
and so
vanishes at $m=0$ for all $b\geqslant 2$. This leaves two values $b=0$ and $b=1$ in the sum,
thus showing that  $\f(\Phi_k)=\f(\Phi^{(0)}_k)$ is given by the required formula.
\epf

\bre\label{rem:typec}\ \ (i)\quad
Due to \cite[Theorem~8.1.5]{f:lc}, we have an isomorphism \eqref{hchiaff}
between the Feigin--Frenkel centre $\z(\wh\g)$ associated with a simple Lie algebra $\g$
and the classical $\Wc$-algebra $\Wc({}^L\g)$.
A direct proof of the isomorphism for the Lie algebras of types $A$, $B$ and $D$ was given in
\cite[Sec.~13.1]{m:so}.
Since the above proof of Theorem~\ref{thm:sphch} does not rely
on the results of \cite{f:lc}, we thus get an independent proof
for type $C$ as well, by
extending
the arguments of \cite[Sec.~13.1]{m:so}. For $\g=\spa_{2n}$ we
get an isomorphism
\ben
\z(\wh\spa_{2n})\cong \Wc(\oa^{}_{2n+1})
\een
by identifying the Harish-Chandra images $\f(\Phi_k)$ with generators of
$\Wc(\oa^{}_{2n+1})$.
\par
(ii)\quad Yet another proof of Theorem~\ref{thm:sphch}
is obtained by starting with
the argument of \cite[Prop.~13.1.14]{m:so}
in the case $k\leqslant n$. The symmetrizer
$S^{(k)}$ is well-defined for these values and the element $(-1)^k\ts \Phi_k$
coincides with the Segal--Sugawara vector $\phi_{kk}$ of \cite{m:ss}
and \cite[Sec.~8.3]{m:so}. Then extend the formula to the region $n+1\leqslant k\leqslant 2n+1$, where
the symmetrizer $S^{(k)}$ may have singularities and cannot be used in the calculation.
Similar to the extrapolation argument of
Sec.~\ref{subsec:ea}, fix a value of $k$ and let $n$ vary, assuming the natural
embeddings of the symplectic Lie algebras
\beql{chainsp}
\spa_{2}\subset \spa_{4}\subset \dots\subset \spa_{2n}\subset \spa_{2n+2}\subset\dots
\eeq
instead of \eqref{chainosp}.
Lemmas~\ref{lem:redu}--\ref{lem:trans}
are valid in the same form in the case $m=0$ so that
the coefficient of any fixed monomial $\mu_{j_1}[r_1]\dots \mu_{j_s}[r_s]$
is a polynomial in $n$ whose degree does not exceed $k$. By the
first part of the proof,
these coefficients are known for all $n\geqslant k$, and hence the
coefficient of the monomial is given by the same formula for all values of $n\geqslant 1$.
\qed
\ere

As shown in \cite[Cor.~13.4.6]{m:so}, a version
of Theorem~\ref{thm:hchfd} for $M=0$,
providing the Harish-Chandra images of the central elements \eqref{ck} in $\U(\spa_{2n})$,
is a straightforward consequence of Theorem~\ref{thm:sphch}.

\appendix

\section{Harish-Chandra images for $\gl_{m|n}$}
\label{sec:gl}

Segal--Sugawara vectors for the general linear Lie superalgebra $\gl_{m|n}$ were constructed
in \cite{mr:mm}, and they were used to describe singular vectors in the Verma
modules at the critical level. The purpose of this appendix is to review these
constructions and calculate the Harish-Chandra images in an explicit form
analogous to Theorem~\ref{thm:hch}.
We will make a connection to the results of the recent paper \cite{afn:ca}
by producing pseudo-differential formulas for such images.
Furthermore, similar
to the orthosymplectic case as in Sec.~\ref{sec:so}, we produce the Harish-Chandra images
of the associated Sugawara operators and construct elements of the quantum Mishchenko--Fomenko
subalgebra of $\U(\gl_{m|n})$; cf. Sec.~\ref{subsec:qmf}.

\subsection{Segal--Sugawara vectors}
\label{subsec:ssvgln}

A standard basis of the general linear Lie superalgebra $\gl_{m|n}$ is formed by elements $E_{ij}$
of the parity $\bi+\bj\mod 2$ for $1\leqslant i,j\leqslant m+n$, with the commutation relations
\eqref{glndr}, where we adapt the notation by setting $\bi=0$ for $1\leqslant i\leqslant m$ and
$\bi=1$ for $m+1\leqslant i\leqslant m+n$.
Now taking $\g=\gl_{m|n}$, consider the
affine Kac--Moody superalgebra \eqref{km}
with the commutation relations
\begin{align}
\non
\big[E_{ij}[r],E_{kl}[s\tss]\tss\big]
=\de_{kj}\ts E_{i\tss l}[r+s\tss]
{}&-\de_{i\tss l}\ts E_{kj}[r+s\tss](-1)^{(\bi+\bj)(\bk+\bl)}\\
\label{commrel}
{}&+K\Big((n-m)\tss\de_{kj}\tss\de_{i\tss l}(-1)^{\bi}
+\de_{ij}\tss\de_{kl}(-1)^{\bi+\bk}\Big)
\ts r\tss\de_{r,-s},
\end{align}
where the element $K$ is even and central, and
we set $E_{ij}[r]=E_{ij}t^r$. The $\ZZ_2$-degree (or parity) of the element $E_{ij}[r]$
is $\bi+\bj\mod 2$.

The {\em vacuum module}
$V_{\cri}(\g)$ at the critical level over $\wh\g$ is defined as the quotient
of the universal enveloping algebra
$\U(\wh\g)$ by the left ideal generated by
$\g[t]$ and $K-1$. It possesses a vertex algebra
structure; see \cite{f:lc} and \cite{k:va}.
The centre $\z(\wh\g)$ of the vertex algebra $V_{\cri}(\g)$
is defined by \eqref{centsup}.
Elements of $\z(\wh\g)$
are called {\em Segal--Sugawara vectors\/}.
The centre
is a commutative associative superalgebra and it can be
identified with
a commutative subalgebra of $\U(t^{-1}\g[t^{-1}])$.

As with the orthosymplectic Lie superalgebras, we consider
the extended Lie superalgebra  $\wh\g\oplus \CC \tau$ with $\tau=-d/dt$.
To reproduce the Segal--Sugawara vectors from \cite{mr:mm},
consider the tensor product superalgebra
\beql{tenpr}
\underbrace{\End\CC^{m|n}\ot\dots\ot\End\CC^{m|n}}_k{}\otimes \U,
\eeq
where $\U=\U(\wh\g\oplus\CC\tau)$.
The symmetric group $\Sym_k$
acts naturally on the tensor product
space $(\CC^{m|n})^{\ot k}$. Along with the symmetrizer defined in \eqref{ha},
introduce the antisymmetrizer
\beql{ak}
a^{(k)}=\frac{1}{k!}\sum_{\si\in\Sym_k}\sgn\si\cdot\si\in\CC\Sym_k
\eeq
and let $H^{(k)}$ and $A^{(k)}$ denote
the respective images of $h^{(k)}$ and $a^{(k)}$
in \eqref{tenpr}. The supertrace is a linear map defined by
\ben
\str:\End \CC^{m|n}\to\CC,\qquad e_{ij}\mapsto \de_{ij}(-1)^{\bi},
\een
where we keep the notation $e_{ij}$ for matrix units in $\End \CC^{m|n}$.

As in \eqref{fra},
for $1 \leqslant a \leqslant k$ and $r \in \ZZ$, define the elements
\beql{efra}
\hE[r]_a \coloneqq \sum_{i,j = 1}^{m+n} 1 ^{\otimes (a-1)} \otimes e_{ij} \otimes 1^{\otimes (k-a)} \otimes \hE_{ij}[r] (-1)^{\bi \bj+\bj},
\eeq
where $\hE_{ij}[r]=\hE_{ij}t^r=E_{ij}(-1)^{\bi}\ts t^r$.
By the results of \cite[Sec.~3.1]{mr:mm},
all the coefficients
$\phi_{k\tss l}, \psi_{k\tss l},\ta_{k\tss l}\in\U(t^{-1}\g[t^{-1}])$
in the expansions
\begin{align}
\str^{}_{1,\dots,k}\ts A^{(k)} \big(\tau+\hE[-1]_1\big)\dots \big(\tau+\hE[-1]_k\big)&={}
\phi_{k\tss 0}\tau^k+\phi_{k1}\tau^{k-1}+\dots+\phi_{kk},
\label{phi}\\[0.4em]
\str^{}_{1,\dots,k}\ts H^{(k)} \big(\tau+\hE[-1]_1\big)\dots \big(\tau+\hE[-1]_k\big)&={}
\psi_{k\tss 0}\tau^k+\psi_{k1}\tau^{k-1}+\dots+\psi_{kk},
\label{psi}\\[0.4em]
\str\ts \big(\tau+\hE[-1]\big)^k&={}
\ta_{k\tss 0}\ts\tau^{k}+\ta_{k1}\ts\tau^{k-1}
+\dots+\ta_{kk},
\label{theta}
\end{align}
are Segal--Sugawara vectors. Moreover, these coefficients can be expressed in terms of the
noncommutative Berezinian through the identities of \cite[Sec.~2.3]{mr:mm}.
It was conjectured in \cite[Remark~3.4(ii)]{mr:mm} that
each of the families $\phi_{kk}$, $\psi_{kk}$ and $\ta_{kk}$
generates the differential algebra $\z(\wh\g)$. A proof of the conjecture in the case $n=1$
was given in \cite{afn:ca}.

\subsection{Harish-Chandra homomorphism}
\label{subsec:hchgln}

Now fix the triangular decomposition \eqref{tridec}
of the Lie superalgebra $\g=\gl_{m|n}$, where the subalgebras $\n_-$ and $\n_+$
are spanned by the elements $E_{ij}$ with $i>j$ and $i<j$, respectively, whereas
$\h$ is spanned by the elements $E_{ii}$ for $i=1,\dots,m+n$.
The {\em affine Harish-Chandra homomorphism}
\beql{hchaffgln}
\f:\U\big(t^{-1}\g[t^{-1}]\big)^{\h}\to \U\big(t^{-1}\h[t^{-1}]\big)
\eeq
is defined in the same way as for the orthosymplectic case in \eqref{hchaff}.
Set
\beql{muire}
\mu_i[-r]=E_{ii}[-r],\qquad i=1,\dots,m+n\fand r=1,2,\dots.
\eeq

We will use the noncommutative elementary and complete supersymmetric functions
in the set of variables $x=(x_1,\dots,x_m,x_{m+1},\dots,x_{m+n})$,
which are defined by the respective formulas
\ben
\bal
e_k(x)&=\sum_{r=0}^k\
\sum_{i_1\geqslant\dots\geqslant i_r\geqslant m+1>i_{r+1}>\dots> i_k}
x_{i_1}\dots x_{i_k},\\[0.3em]
h_k(x)&=\sum_{r=0}^k\
\sum_{i_1\leqslant\dots\leqslant i_r\leqslant
m<i_{r+1}<\dots< i_k}
x_{i_1}\dots x_{i_k}.
\eal
\een
We also set
\ben
p_k(x)=\sum_{r=1}^k(-1)^{r-1}\ts r\ts h_{k-r}(x)\ts e_r(x).
\een

\bth\label{thm:hcha}
For all $k\geqslant 1$,
the Harish-Chandra images of the polynomials \eqref{phi},
\eqref{psi} and \eqref{theta} are given by the respective formulas
\begin{align}
{}&e_k\big(\tau+\mu_{1}[-1],\dots,\tau+\mu_{m}[-1],-\tau+\mu_{m+1}[-1],\dots,-\tau+\mu_{m+n}[-1]\big),
\label{hchphi}\\[0.3em]
{}&h_k\big(\tau+\mu_{1}[-1],\dots,\tau+\mu_{m}[-1],-\tau+\mu_{m+1}[-1],\dots,-\tau+\mu_{m+n}[-1]\big),
\label{hchpsi}
\\[0.3em]
\label{hchth}
{}&p_k\big(\tau+\mu_{1}[-1],\dots,\tau+\mu_{m}[-1],-\tau+\mu_{m+1}[-1],\dots,-\tau+\mu_{m+n}[-1]\big).
\end{align}
\eth

\bpf
The arguments are parallel to the respective proofs of \cite[Thm~3.10, Cor.~3.11 and Cor.~3.12]{mr:mm}.
Formula \eqref{hchphi} is immediate from \cite[Prop.~2.3]{mr:mm}.
The Berezinian identities of \cite[Thm~2.13]{mr:mm} imply
\begin{align}\label{charfermthm}
\Ber(1+q\tss T)&=\sum_{k=0}^{\infty} q^k\ts
\str^{}_{1,\dots,k}\ts A^{(k)} T_1\dots T_k,\\
\label{charbosthm}
\big[\Ber(1-q\tss T)\big]^{-1}&=\sum_{k=0}^{\infty}
q^k\ts\str^{}_{1,\dots,k}\ts H^{(k)} T_1\dots T_k,\\
\label{newtonthm}
\big[\Ber(1+q\tss T)\big]^{-1}\ts\di_q\tss \Ber(1+q\tss T)
&=\sum_{k=0}^{\infty} (-q)^k\ts\str\ts T^{k+1},
\end{align}
where we set $T=\tau+\hE[-1]$, and $q$ is a variable. By \eqref{hchphi} we have
\begin{multline}
\f:\Ber(1+q\tss T)
\mapsto \big(1+q(\tau-\mu_{m+n}[-1])\big)^{-1}\dots \big(1+q(\tau-\mu_{m+1}[-1])\big)^{-1}\\[0.4em]
{}\times \big(1+q(\tau+\mu_{m}[-1])\big)\dots \big(1+q(\tau+\mu_{1}[-1])\big).
\non
\end{multline}
Now \eqref{hchpsi} and \eqref{hchth} follow from \eqref{charbosthm} and \eqref{newtonthm},
respectively.
\epf

\subsection{Pseudo-differential operators}
\label{subsec:pd}

Following \cite{afn:ca}, we will obtain pseudo-differential operator formulas
for the Harish-Chandra images of Theorem~\ref{thm:hcha}.
Consider the constant term Segal--Sugawara vectors in \eqref{phi} and \eqref{psi}
and set
\ben
\Ec_k=\f(\phi_{kk})\Fand \Hc_k=\f(\psi_{kk}).
\een
These are the respective constant terms of the polynomials in $\tau$
appearing in \eqref{hchphi} and \eqref{hchpsi}.

\bco\label{cor:pseudo}
We have the pseudo-differential operator expansions
\ben
\bal
\big(\tau-\mu_{m+n}[-1]\big)^{-1}\dots &\big(\tau-\mu_{m+1}[-1]\big)^{-1}\\
\times{\ts}&\big(\tau+\mu_{m}[-1]\big)\dots \big(\tau+\mu_{1}[-1]\big)=\sum_{r=0}^{\infty}\Ec_r \tau^{m-n-r},
\eal
\een
and
\ben
\bal
\big(\tau+\mu_{1}[-1]\big)^{-1}\dots &\big(\tau+\mu_{m}[-1]\big)^{-1}\\
\times{\ts}&\big(\tau-\mu_{m+1}[-1]\big)\dots \big(\tau-\mu_{m+n}[-1]\big)&=\sum_{r=0}^{\infty}(-1)^r\Hc_r \tau^{-m+n-r}.
\eal
\een
\eco

\bpf
Note that the relations
\beql{taumu}
\big[\tau,\mu_i[r]\tss\big]=-r\ts \mu_i[r-1],
\eeq
will remain valid under the replacement $\tau\mapsto q^{-1}+\tau$.
Therefore, is is enough to prove the expansions obtained after this
replacement.
The first pseudo-differential operator then takes the form
\begin{multline}
q^{-m+n}\tss\big(1+q(\tau-\mu_{m+n}[-1])\big)^{-1}\dots \big(1+q(\tau-\mu_{m+1}[-1])\big)^{-1}\\[0.4em]
{}\times \big(1+q(\tau+\mu_{m}[-1])\big)\dots \big(1+q(\tau+\mu_{1}[-1])\big).
\non
\end{multline}
As we pointed out in the proof of Theorem~\ref{thm:hcha}, this
equals $q^{-m+n}\ts\f(\Ber(1+q\tss T))$ and so coincides with
the generating function of the polynomials \eqref{hchphi} multiplied by $q^{-m+n}$.
On the other hand, we have the expansion
\beql{ekek}
\bal
e_k\big(\tau+\mu_{1}[-1],\dots,\tau+\mu_{m}[-1],-\tau+\mu_{m+1}[-1],&\dots,-\tau+\mu_{m+n}[-1]\big)\\
{}&=\sum_{r=0}^k\binom{m-n-r}{k-r}\tss\Ec_r\tss\tau^{k-r},
\eal
\eeq
which is verified in the same way as in the cases $m=0$ and $n=0$; cf. \cite[Prop.~12.4.4 and Prop.~12.4.7]{m:so}. Hence, the pseudo-differential operator becomes
\ben
\sum_{k=0}^{\infty}\ts \sum_{r=0}^k
\binom{m-n-r}{k-r}\tss\Ec_r\tss\tau^{k-r}\ts q^{-m+n+k}
=\sum_{r=0}^{\infty}\ts\Ec_r\ts\sum_{k=r}^{\infty}\binom{m-n-r}{k-r}q^{-m+n+k}\tau^{k-r}
\een
which coincides with
\ben
\sum_{r=0}^{\infty}\Ec_r\ts (q^{-1}+\tau)^{m-n-r},
\een
thus proving the first expansion. To verify the second,
we will use the replacement $\tau\mapsto -q^{-1}+\tau$ and argue in a similar way.
The second pseudo-differential operator takes the form
\begin{multline}
(-q)^{m-n}\tss\big(1-q(\tau+\mu_{1}[-1])\big)^{-1}\dots \big(1-q(\tau+\mu_{m}[-1])\big)^{-1}\\[0.4em]
{}\times \big(1+q(-\tau+\mu_{m+1}[-1])\big)\dots \big(1+q(-\tau+\mu_{m+n}[-1])\big)
\non
\end{multline}
which coincides with the Harish-Chandra image of the series $\big[\Ber(1-q\tss T)\big]^{-1}$
multiplied by $(-q)^{m-n}$. Similar to \eqref{ekek}, we have
\ben
\bal
h_k\big(\tau+\mu_{1}[-1],\dots,\tau+\mu_{m}[-1],-\tau+\mu_{m+1}[-1],&\dots,-\tau+\mu_{m+n}[-1]\big)\\
{}&=\sum_{r=0}^k\binom{m-n+k-1}{k-r}\tss\Hc_r\tss\tau^{k-r}.
\eal
\een
Hence, the second pseudo-differential operator becomes
\ben
(-1)^{m-n}\ts\sum_{k=0}^{\infty}\ts \sum_{r=0}^k
\binom{m-n+k-1}{k-r}\tss\Hc_r\tss\tau^{k-r}\ts q^{m-n+k}
\een
which equals
\ben
(-1)^{m-n}\sum_{r=0}^{\infty}\ts\Hc_r\ts\sum_{k=r}^{\infty}\binom{m-n+k-1}{k-r}q^{m-n+k}\tau^{k-r}
=\sum_{r=0}^{\infty}(-1)^r\ts\Hc_r\ts (-q^{-1}+\tau)^{-m+n-r},
\een
thus completing the proof.
\epf

\bre\label{rem:afnpaper}
The pseudo-differential operator expansions were used in \cite{afn:ca} (in the case $n=1$)
in a slightly different settings. In our notation, for arbitrary $m$ and $n$,
define polynomials
$W_k$
in the variables $\mu_i[-r]$ by the expansion
of the pseudo-differential operator
\ben
\big(\tau-\mu_1[-1]\big)\dots \big(\tau-\mu_m[-1]\big)
\big(\tau+\mu_{m+1}[-1]\big)^{-1}\dots \big(\tau+\mu_{m+n}[-1]\big)^{-1}
=\sum_{k=0}^{\infty} W_k\ts \tau^{m-n-k}.
\een
As in the proof of Corollary~\ref{cor:pseudo},
we can derive
the following explicit expressions for the coefficients $W_k$:
\ben
W_k={\bar e}_k\big(\tau-\mu_{1}[-1],\dots,\tau-\mu_{m}[-1],-\tau-\mu_{m+1}[-1],
\dots,-\tau-\mu_{m+n}[-1]\big)\ts 1,
\een
assuming that $\tau\ts 1=0$, where a different version
of the noncommutative elementary supersymmetric functions is used, namely
\ben
{\bar e}_k(x)=\sum_{r=0}^k\
\sum_{i_1<\dots< i_r<
m+1\leqslant i_{r+1}\leqslant\dots\leqslant i_k}
x_{i_1}\dots x_{i_k}.
\een
Therefore, by \cite[Prop.~2.3 and Remark~2.4(i)]{mr:mm}, the coefficient $W_k$
coincides with the Harish-Chandra image of the Segal--Sugawara vector $\phi_{kk}$
under a modified version of the homomorphism \eqref{hchaffgln}. Namely, the roles
of the subalgebras $\n_+$ and $\n_-$ should be swapped, together with the
sign change $\mu_i[-r]\mapsto -\mu_i[-r]$.
\qed
\ere

\subsection{Sugawara operators}
\label{sec:sogln}

The definition \eqref{complua} of the completed universal enveloping algebra $\wt\U_{\cri}(\wh\g)$
extends to the case $\g=\gl_{m|n}$ without any essential changes, except that the critical level
is now $K=1$.
We will use Theorem~\ref{thm:hcha} to calculate the
Harish-Chandra images of the associated
{\em Sugawara operators} which we regard as
elements of the centre $\Z(\wh\g)$ of $\wt\U_{\cri}(\wh\g)$; cf. Corollary~\ref{cor:hchss}.

For any $i,j\in\{1,\dots,m+n\}$, introduce the Laurent series $E_{ij}(u)$
with coefficients in $\U_{\cri}(\wh\g)$ by
\ben
E_{ij}(u)=\sum_{r\in\ZZ} E_{ij}[r]\tss u^{-r-1}
\een
and combine them into the matrix $\hE(u)=[\hE_{ij}(u)]$ with $\hE_{ij}(u)=E_{ij}(u)(-1)^{\bi}$.
Set $T(u)=\di_u+\wh E(u)$, where $\di_u$ is understood
as a scalar matrix of size $m+n$.
The following proposition is
a reformulation of \cite[Cor.~3.5]{mr:mm}, where we use the normally ordered products
as recalled in Section~\ref{sec:so}.

\bpr\label{prop:sugagln}
All coefficients of the Laurent series
$\phi_{kl}(u)$, $\psi_{kl}(u)$ and $\ta_{kl}(u)$
defined by the decompositions
\begin{align}
:\str^{}_{1,\dots,k}\ts A^{(k)} T(u)_1\dots T(u)_k:{}&=
\phi_{k\tss 0}(u)\ts\di_u^{\tss k}+\phi_{k1}(u)\ts\di_u^{\tss k-1}
+\dots+\phi_{kk}(u),\label{phiso}\\[0.4em]
:\str^{}_{1,\dots,k}\ts H^{(k)} T(u)_1\dots T(u)_k:{}&=
\psi_{k\tss 0}(u)\ts\di_u^{\tss k}+\psi_{k1}(u)\ts\di_u^{\tss k-1}
+\dots+\psi_{kk}(u),\label{psiso}\\[0.4em]
:\str\ts T(u)^k:{}&
= \ta_{k\tss 0}(u)\ts\di_u^{\tss k}+\ta_{k1}(u)\ts\di_u^{\tss k-1} +\dots
+\ta_{kk}(u),\label{taso}
\end{align}
are Sugawara operators for $\gl_{m|n}$.
\qed
\epr

The definition \eqref{hchcompl} of the affine Harish-Chandra homomorphism
extends to the case $\g=\gl_{m|n}$ providing a homomorphism
\ben
\f:\Z(\wh\g)\to \wt\Pi
\een
from the centre $\Z(\wh\g)$ of the completed universal enveloping algebra $\wt\U_{\cri}(\wh\g)$
to the completion $\wt\Pi$ of the symmetric algebra $\Pi=\Sr\big(\h\tss[t,t^{-1}]\big)$.
Introduce notation for the basis elements of $\h\tss[t,t^{-1}]$ by setting
\beql{muirfullgln}
\mu_i[r]=E_{ii}[r]\Fand \mu_i(u)=\sum_{r\in\ZZ} \mu_i[r]\tss u^{-r-1},
\eeq
where $i=1,\dots,m+n$.
The following corollary is immediate from Theorem~\ref{thm:hcha}.

\bco\label{cor:hchssgln}
For all $k\geqslant 1$,
the Harish-Chandra images of the polynomials \eqref{phiso},
\eqref{psiso} and \eqref{taso} are given by the respective formulas
\begin{align}
{}&e_k\big(\di_u+\mu_{1}(u),\dots,\di_u+\mu_{m}(u),-\di_u+\mu_{m+1}(u),\dots,-\di_u+\mu_{m+n}(u)\big),
\non\\[0.3em]
{}&h_k\big(\di_u+\mu_{1}(u),\dots,\di_u+\mu_{m}(u),-\di_u+\mu_{m+1}(u),\dots,-\di_u+\mu_{m+n}(u)\big),
\non\\[0.3em]
{}&p_k\big(\di_u+\mu_{1}(u),\dots,\di_u+\mu_{m}(u),-\di_u+\mu_{m+1}(u),\dots,-\di_u+\mu_{m+n}(u)\big).
\non
\end{align}
\eco

\subsection{Quantum Mishchenko--Fomenko subalgebras}
\label{subsec:qmfa}

As with the orthosymplectic Lie superalgebras discussed in Section~\ref{subsec:qmf},
suppose that $P\in\Sr(\g)^{\g}$ for $\g=\gl_{m|n}$
is a $\g$-invariant of $\Sr(\g)$ under the adjoint action.
Take any element $\mu\in\g^*$ which
vanishes on the odd elements of $\g$. Regarding $P$
as a polynomial in the $E_{ij}$, use a `shift of argument'
to replace the variables by
$
E_{ij}\mapsto E_{ij}+t\tss \mu(E_{ij}),
$
where $t$ is a variable. After this replacement, the new polynomial expands
as a polynomial in $t$,
\beql{poltea}
P^{}_{(0)}+P^{}_{(1)}\tss t+\dots+P^{}_{(k)}\tss t^k,
\eeq
thus defining
elements $P^{}_{(i)}\in \Sr(\g)$ associated with $P$ and $\mu$.
The {\em (classical) Mishchenko--Fomenko subalgebra}
$\overline\Ac_{\mu}$ of $\Sr(\g)$
is generated by all elements $P^{}_{(i)}$
associated with all $\g$-invariants $P\in \Sr(\g)^{\g}$.
The subalgebra $\overline\Ac_{\mu}$ is Poisson super-commutative
with respect to the Lie--Poisson super-bracket on $\Sr(\g)$; cf. \cite{r:si}
and \cite[Sec.~9.1]{m:so}.
{\em Vinberg's quantization problem}
\cite{v:sc}
asks whether it is possible to
construct a commutative
subalgebra $\Ac_{\mu}$ of $\U(\g)$ such that $\gr\Ac_{\mu}=\overline\Ac_{\mu}$.

Define the subalgebra $\Ac_{\mu}\subset\U(\g)$ as the image
of $\z(\wh\g)$ with respect to the
homomorphism
\beql{evalrgmoda}
\U\big(t^{-1}\g[t^{-1}]\big)\to \U(\g),
\qquad E_{ij}[-r]\mapsto E_{ij}\tss z^{-r}+\de_{r1}\ts\mu(E_{ij}),\quad r>0,
\eeq
where $z\in\CC$ is nonzero.
The image is easily seen not to depend on $z$.
Explicit formulas for elements of $\Ac_{\mu}$ then
follow from the results of \cite{mr:mm}.
Introduce the matrix $\mu=[\mu(\hE_{ij})]$.

\bpr\label{prop:qmfa}
The coefficients of all polynomials in $z^{-1}$ given by
\ben
\bal
\str^{}_{1,\dots,k}\ts A^{(k)} &\big(\di_z+\mu_1+\hE_1\tss z^{-1}\big)
\dots \big(\di_z+\mu_k+\hE_k\tss z^{-1}\big)\ts 1,
\\[0.4em]
\str^{}_{1,\dots,k}\ts H^{(k)} &\big(\di_z+\mu_1+\hE_1\tss z^{-1}\big)
\dots \big(\di_z+\mu_k+\hE_k\tss z^{-1}\big)\ts 1,
\eal
\een
and
\ben
\str\ts \big(\di_z+\mu+\hE\tss z^{-1}\big)^k \ts 1,
\een
with positive integer values of $k$,
belong to the super-commutative superalgebra $\Ac_{\mu}$.
\qed
\epr

If the conjecture of
\cite[Remark~3.4(ii)]{mr:mm}
on generators of $\z(\wh\g)$ holds, then
the elements defined in Proposition~\ref{prop:qmfa} generate the superalgebra
$\Ac_{\mu}$. Since the conjecture is proven in the case
$n=1$ in \cite{afn:ca}, we know this holds for $n=1$.
Moreover, we expect the following to hold.

\bcj\label{conj:mfa}
The superalgebras $\Ac_{\mu}$ solve Vinberg's quantization problem:
$\gr\Ac_{\mu}=\overline\Ac_{\mu}$ for all $\mu$.
\ecj

\section*{Declarations}

\subsection*{Competing interests}
The authors have no competing interests to declare that are relevant to the content of this article.

\subsection*{Acknowledgements}
Our work was supported by the Australian Research Council, grant DP240101572.
E.R. wishes to thank the School of Mathematics and Statistics
at the University of Sydney
for the warm hospitality during his visit.

\subsection*{Availability of data and materials}
No data was used for the research described in the article.


\bigskip
\bigskip

\small
\noindent
School of Mathematics and Statistics\newline
University of Sydney,
NSW 2006, Australia\newline
{\tt alexander.molev@sydney.edu.au}

\vspace{5 mm}

\noindent
School of Mathematics and Statistics\newline
University of Sydney,
NSW 2006, Australia\newline
{\tt madeline.nurcombe@sydney.edu.au}

\vspace{5 mm}

\noindent
Laboratoire d'Annecy de Physique Th\'{e}orique,
CNRS and Universit\'{e} de Savoie\newline
BP 110, 74941 Annecy-le-Vieux Cedex, France\newline
{\tt eric.ragoucy@lapth.cnrs.fr}

\end{document}